\documentclass[12pt, titlepage]{amsart}
\usepackage{amssymb,mathrsfs,bm,setspace}
\usepackage[author-year]{amsrefs}
\title[Local Times and SPDEs]{%
	A Local-Time Correspondence for\\Stochastic 
	Partial Differential Equations}
\thanks{The research of D.K.
	is supported in part by the National Science Foundation
	grant DMS-0404729}
\author[Foondun]{Mohammud Foondun}
	\address{Department of Mathematics, The University of Utah,
		155 S. 1400 E. Salt Lake City, UT 84112--0090, USA}
	\email{mohammud@math.utah.edu}
	\urladdr{http://www.math.utah.edu/\~{}mohammud}
\author[Khoshnevisan]{Davar Khoshnevisan}
	\address{Department of Mathematics, The University of Utah,
		155 S. 1400 E. Salt Lake City, UT 84112--0090, USA}
	\email{davar@math.utah.edu}
	\urladdr{http://www.math.utah.edu/\~{}davar}
\author[Nualart]{Eulalia Nualart}
	\address{Institut Galil\'ee, Universit\'e
		Paris 13, 93430 Villetaneuse, France}
	\email{nualart@math.univ-paris13.fr}
	\urladdr{http://www.math.univ-paris13.fr/\~{}nualart}

\date{November 11, 2007}
\theoremstyle{plain}{
\newtheorem{theorem}{Theorem}[section]}
\theoremstyle{plain}{
\newtheorem{proposition}[theorem]{Proposition}}
\theoremstyle{plain}{
   \newtheorem{lemma}[theorem]{Lemma}}
\theoremstyle{plain}{
   \newtheorem{corollary}[theorem]{Corollary}}
\theoremstyle{definition}{
   \newtheorem{definition}[theorem]{Definition}}
\theoremstyle{definition}{
   \newtheorem{example}[theorem]{Example}}
\theoremstyle{remark}{
   \newtheorem{remark}[theorem]{Remark}}
\theoremstyle{definition}{
   }
\theoremstyle{plain}

\numberwithin{equation}{section}
\newcommand{\R}{\mathbf{R}}
\renewcommand{\P}{\mathrm{P}}
\newcommand{\E}{\mathrm{E}}
\newcommand{\1}{{\bf 1}}
\newcommand{\F}{\mathscr{F}}
\renewcommand{\Re}{\mathrm{Re}\,}

\renewcommand{\L}{{{\mathscr{L}}}}
\newcommand{\sE}{\mathscr{E}}
\newcommand{\sS}{\mathscr{S}}
\renewcommand{\phi}{\varphi}
\begin{document}
\onehalfspacing
\begin{abstract}
	It is frequently the case that a white-noise-driven
	parabolic and/or hyperbolic stochastic partial differential
	equation (SPDE)
	can have random-field solutions only in spatial dimension
	one. Here we show that in many cases, where
	the ``spatial operator'' is the $L^2$-generator of a L\'evy process $X$,
	a linear SPDE has a random-field solution if and only if
	the symmetrization of $X$ possesses local times. 
	This result gives a probabilistic
	reason for the lack of existence of random-field solutions in 
	dimensions strictly bigger than one.

	In addition, we prove that the solution to the SPDE is [H\"older]
	continuous in its spatial variable if and only if the said
	local time is [H\"older] continuous in its spatial variable.
	We also produce examples where the random-field solution exists,
	but is almost surely
	unbounded in every open subset of space-time.
	Our results are based on first establishing
	a quasi-isometry between the
	linear $L^2$-space of the
	weak solutions of a family of linear SPDEs, on one hand,
	and the Dirichlet space generated by the symmetrization
	of $X$, on the other hand.
	
	We study mainly linear equations in order to present the
	local-time correspondence at a modest technical level. However,
	some of our work has consequences for nonlinear SPDEs as well.
	We demonstrate this assertion by studying a family of
	parabolic SPDEs that have additive nonlinearities. 
	For those equations we prove that if
	the linearized problem has a random-field solution,
	then so does the nonlinear SPDE. Moreover, the solution
	to the linearized equation is [H\"older] continuous if and only if
	the solution to the nonlinear equation is. And the solutions
	are bounded and unbounded together as well. Finally, we prove
	that in the cases that the solutions are unbounded, they 
	almost surely blow up at exactly the same points.
\end{abstract}
\maketitle
\tableofcontents

\section{\bf Introduction}

We consider the \emph{stochastic heat equation}
inspired by the fundamental works of 
Pardoux \ycites{Pardoux:75b,Pardoux:75a,Pardoux:72},
Krylov and Rozovskii \ycites{KR:79b,KR:79a,KR:77},
and Funaki \ycite{Funaki}. Let
$\dot{w}$ denote space-time white noise, $t$ is nonnegative,
$x$ is in $\R^d$, and the Laplacian acts on the $x$ variable.
Then, we have:
\begin{equation}\label{heat}\left| \begin{split}
	&\partial_t  H (t\,,x) =
		(\Delta  H )(t\,,x) + \dot{w}(t\,,x),\\
	& H (0\,,x) =0.
\end{split}\right.\end{equation}
Let us consider also the \emph{stochastic wave
equation} of Caba\~{n}a \ycite{Cabana}:
\begin{equation}\label{wave}\left|\begin{split}
	&\partial_{tt}  W (t\,,x) =
		(\Delta  W )(t\,,x) + \dot{w}(t\,,x),\\
	& W (0\,,x) = \partial_t  W (0\,,x) = 0.
\end{split}\right.\end{equation}

One of the common features of \eqref{heat} and \eqref{wave}
is that they suffer from a \emph{curse of dimensionality}.
Namely, these equations can have random-field solutions
only in dimension one. Moreover, this curse of dimensionality
appears to extend beyond the linear parabolic setting of
\eqref{heat}, or the linear hyperbolic setting of
\eqref{wave}. For instance, see
Perkins \ycite{Perkins}*{Corollary III.4.3}
for an example from superprocesses, and 
Walsh \ycite{Walsh}*{Chapter 9} for an example
from statistical mechanics. 

One can informally ascribe this
curse of dimensionality to the ``fact'' that
while the Laplacian smooths, white noise
roughens. In one dimension,
the roughening effect of white noise turns out to be small relative
to the smoothing properties of the Laplacian, and we thus have
a random-field solution. However, in dimensions greater than
one white noise is much too rough, and the Laplacian
cannot smooth the solution enough to yield a
random field.

Dalang and Frangos \ycite{DalangFrangos} were able to
construct a first fully-rigorous explanation of the
curse of dimensionality.
They do so by first replacing white noise by a Gaussian noise that
is white in time and colored in space. And then they describe
precisely the roughening effect of the noise on the solution,
viewed as a random generalized function. See also 
Brze{\'z}niak and van Neerven \ycite{BvN},
Dalang and Mueller \ycite{DalangMueller},
Millet and Sanz-Sol\'e \ycite{MilletSanzsole},
Peszat \ycite{Peszat}, and
Peszat and Zabczyk \ycite{PeszatZabczyk}.
More recently,
Dalang and Sanz-Sol\'e \ycite{DalangSole:05}
study fully nonlinear stochastic wave equations driven
by noises that are white in 
time and colored in space, and operators that are arbitrary 
powers of the Laplacian.

In this article we present a different explanation of this phenomenon.
Our approach is to describe accurately
the smoothing effect of the Laplacian
in the presence of white noise. 
Whereas the answer of
Dalang and Frangos \ycite{DalangFrangos}
is analytic, ours is  probabilistic.
For instance, we will see soon that \eqref{heat} and \eqref{wave}
have solutions only in dimension one because $d$-dimensional Brownian motion
has local times only in dimension one
[Theorem \ref{th:exist}]. Similarly, when $d=1$,
the solution to \eqref{heat} [and/or \eqref{wave}] is continuous
in $x$ because the local time of one-dimensional Brownian motion
is continuous in its spatial variable.

The methods that we employ also
give us a local-time paradigm that makes precise the claim that
the stochastic PDEs \eqref{heat}
and \eqref{wave} ``have random-field solutions in dimension
$d=2-\epsilon$ for all $\epsilon\in(0\,,2)$.''
See Example \ref{example:fractal} below, where we introduce a
family of SPDEs on fractals.

An outline of this paper follows:
In \S\ref{sec2} we describe a suitable generalization of
the stochastic PDEs \eqref{heat} and \eqref{wave} that have
sharp local-time correspondences. This section contains the main
results of the paper. The existence theorem of
\S\ref{sec2} is proved in
\S\ref{sec3}. 
Section \ref{sec4} contains the proof of
the necessary and sufficient local-time condition for continuity
of our SPDEs in their space variables. 
Although we are not aware
of any interesting connections between local times and 
temporal regularity of the solutions
of SPDEs, we have included \S\ref{sec5}
that contains a sharp analytic condition
for temporal continuity of the SPDEs of \S\ref{sec2}.
We also produce examples
of SPDEs that have random-field solutions which are
almost surely unbounded in every open space-time set
[Example \ref{example:sharp}].

Section \ref{sec6} discusses issues
of H\"older continuity in either space, or time, variable.
In  \S\ref{sec7} we establish a very general (but somewhat
weak) connection between Markov processes, their local times,
and solutions to various linear SPDEs. The material of that section
is strongly motivated by the recent article of
Da Prato \ycite{DaPrato} who
studies Kolmogorov SPDEs that are not unlike those studied here, but
also have multiplicative nonlinearities.
We go on to produce examples
where one can make sense of a random-field solution 
to \eqref{heat} in dimension $2-\epsilon$ for all $\epsilon\in(0\,,2)$.
These solutions in fact turn out to be jointly H\"older
continuous, but we will not dwell on that here.
Finally, in \S\ref{sec8}
we discuss a parabolic
version of the SPDEs of \S\ref{sec2} that have additive
nonlinearities;
we prove the existence of solutions and describe exactly when and where
these solutions blow up. 

The chief aim of this paper is to point out various
interesting and deep connections between the local-time theory of
Markov processes and families of stochastic partial differential
equations. In all cases, we have strived to study the simplest SPDEs
that best highlight these connections. But it would also
be interesting to study much more general equations.

There are other connections between
local times of Markov processes and Gaussian processes that
appear to be different from those presented here. 
For a sampler of those isomorphism theorems see
Brydges, Fr\"ohlich, and Spencer \ycite{BFS},
Dynkin \ycite{Dynkin}, and 
Eisenbaum \ycite{Eisenbaum}.
Marcus and Rosen \ycite{MarcusRosen} contains an
excellent and complete account
of this theory. Diaconis and Evans \ycite{DiaconisEvans}
have introduced yet a different isomorphism theorem. 

Finally, we conclude by mentioning what we mean by
``local times,''  as there are many [slightly] different
versions in the literature. Given a stochastic
process $Y:=\{Y_t\}_{t\ge 0}$ on $\R^d$, consider the
occupation measure[s],
\begin{equation}\label{OM}
	Z(t\,,\phi) := \int_0^t \phi(Y_s)\, ds
	\qquad\text{for all $t\ge 0$ and 
	measurable $\phi:\R^d\to\R_+$.}
\end{equation}
We can identify each $Z(t\,,\bullet)$ with a measure
in the usual way.
Then, we say that $Y$ has local times when $Z(t\,,dx)\ll
dx$ for all $t$. The local times of $Y$ are themselves
defined by $Z(t\,,x) := Z(t\,,dx)/dx$. It follows that if
$Y$ has local times, then $Z(t\,,\phi) = 
\int_{\R^d} Z(t\,,x) \phi(x)\, dx$ a.s.\
for every $t\ge 0$ and all
measurable functions $\phi:\R^d\to\R_+$.
And the converse holds also.

\section{\bf More general equations}\label{sec2}
In order to describe when \eqref{heat}
and \eqref{wave} have random-field solutions, and why, we
study more general equations. 

Let $\L $ denote the generator of a $d$-dimensional
L\'evy process $X:=\{X_t\}_{t\ge 0}$
with characteristic exponent $\Psi$. We can normalize things
so that
$\E\exp(i\xi\cdot X_t)=\exp(-t\Psi(\xi))$, and consider
$\L $ as an $L^2$-generator with domain
\begin{equation}
	\text{Dom}\,\L  := \left\{ f\in L^2(\R^d):\
	\int_{\R^d} |\hat f(\xi)|^2\, |\Psi(\xi)|\, d\xi<\infty\right\}.
\end{equation}
As usual, $\hat{f}$ denotes the Fourier transform of
$f$; we opt for the normalization
\begin{equation}
	\hat{f}(\xi):=
	\int_{\R^d} e^{i\xi\cdot x} f(x)\, dx
	\qquad\text{for all $f\in L^1(\R^d)$.}
\end{equation}
[The general $L^2$-theory of Markov processes
is described in great depth in
Fukushima, \={O}shima, and Takeda \ycite{FOT} in the symmetric
case, and Ma and R\"ockner \ycite{MaRockner} for the general case.]

In this way, we can---and will---view $\L $ as a generalized
convolution operator with Fourier multiplier $\hat{\L }(\xi):=
-\overline{\Psi(\xi)}$. 

We consider two families of stochastic partial differential
equations. The first is the stochastic heat equation
for $\L $:
\begin{equation}\label{heat:L}\left|\begin{split}
	&\partial_t  H (t\,,x) =
		(\L   H )(t\,,x) + \dot{w}(t\,,x),\\
	& H (0\,,x) =0,
\end{split}\right.\end{equation}
where $x$ ranges over $\R^d$ and $t$ over 
$\R_+:=[0\,,\infty)$,
and the operator $\L $ acts on the $x$ variable.

We also consider hyperbolic SPDEs of the following wave type:
\begin{equation}\label{wave:L}\left|\begin{split}
	&\partial_{tt}  W (t\,,x) =
		(\L   W )(t\,,x) + \dot{w}(t\,,x),\\
	& W (0\,,x) =\partial_t  W (0\,,x)=0.
\end{split}\right.\end{equation}

Recall that $X:=\{X_t\}_{t\ge 0}$ is a L\'evy process whose generator
is $\L $. Let $X'$ denote an independent copy
of $X$ and define its symmetrization, a la Paul L\'evy, by
\begin{equation}
	\bar{X}_t := X_t - X_t'\qquad\text{for all
	$t\ge 0$}.
\end{equation}
It is a standard fact that $\bar{X}$ is a symmetric
L\'evy process with characteristic exponent $2\Re\Psi$.

The following is our first main result, though the
precise meaning of its terminology is yet to be explained.

\begin{theorem}\label{th:exist}
	The stochastic heat equation \eqref{heat:L} has
	random-field solutions if and only if the
	symmetric L\'evy process $\bar{X}$ has local times.
	The same is true for the stochastic wave equation
	\eqref{wave:L}, provided that the process $X$
	is itself symmetric.
\end{theorem}

\begin{remark}
	Brze{\'z}niak and van Neerven \ycite{BvN}
	consider parabolic equations of the type
	\eqref{heat:L}, where $\L$ is a pseudo-differential
	operator with a symbol that is bounded below
	in our language [bounded above in theirs]. If we apply
	their theory with our constant-symbol
	operator $\L$, then their bounded-below condition
	is equivalent to the symmetry of the L\'evy process
	$X$, and their main result is equivalent to the
	parabolic portion of Theorem \ref{th:exist},
	under the added condition that $X$ is symmetric.
	\qed
\end{remark}

It is well-known that when $d\ge 2$,
L\'evy processes in $\R^d$ \emph{do not}
have local times \cite{Hawkes}. It follows from this that 
neither \eqref{heat:L} nor \eqref{wave:L}
[under a symmetry assumption on $X$] can ever have
random-field solutions in dimension greater than one.
But it is possible that there are no random-field
solutions even in dimension one. Here is one such example;
many others exist.
\begin{example}\label{example:exist:SS}
	Suppose $X$ is a strictly stable process in $\R$
	with stability index $\alpha\in(0\,,2]$.
	It is very well known that 
	the symmetric L\'evy process $\bar{X}$ has local times
	if and only if $\alpha>1$; see, for example,
	Hawkes \ycite{Hawkes}. If $X$ is itself
	symmetric, then $\L =-(-\Delta)^{\alpha/2}$
	is the $\alpha$-dimensional fractional Laplacian
	\cite{Stein}*{Chapter V, \S1.1}, and also the stochastic
	wave equation \eqref{wave:L}
	has a random-field solution if and
	only if $\alpha> 1$.\qed
\end{example}

The local-time correspondence of Theorem \ref{th:exist}
is not a mere accident.
In fact, the next two theorems suggest far deeper connections
between the solutions to the linear SPDEs of this paper
and the theory of local times of Markov processes. 
We emphasize that the next 
two theorems assume
the existence of a random-field solution to one
of the stochastic PDEs \eqref{heat:L} and/or \eqref{wave:L}.
Therefore, they are inherently one dimensional statements.

\begin{theorem}\label{th:cont}
	Assume that $d=1$ and
	the stochastic heat equation \eqref{heat:L} has
	a random-field solution $\{  H (t\,,x)\}_{t\ge 0,x\in\R}$.
	Then, the following are equivalent:
	\begin{enumerate}
		\item There exists $t> 0$ such that $x\mapsto  H (t\,,x)$
			is a.s.\ continuous.
		\item For all $t> 0$, $x\mapsto  H (t\,,x)$ is a.s.\
			continuous.
		\item The local times of $\bar{X}$
			are a.s.\ continuous in their spatial variable.
	\end{enumerate}
	The same equivalence is true for the solution $ W $
	to the stochastic wave equation
	\eqref{wave:L}, provided that the process $X$
	is itself symmetric.
\end{theorem}

\begin{theorem}\label{th:Holder}
	Assume that $d=1$ and
	the stochastic heat equation \eqref{heat:L} has
	a random-field solution $\{  H (t\,,x)\}_{t\ge 0,x\in\R}$.
	Then, the following are equivalent:
	\begin{enumerate}
		\item There exists $t> 0$ such that $x\mapsto  H (t\,,x)$
			is a.s.\ H\"older continuous.
		\item For all $t> 0$, $x\mapsto  H (t\,,x)$ is a.s.\
			H\"older continuous.
		\item The local times of $\bar{X}$
			are a.s.\ H\"older continuous in their spatial variable.
	\end{enumerate}
	If the process $X$ is itself symmetric, then the preceding
	conditions are also equivalent to the H\"older
	continuity of the solution
	$ W $ to the stochastic wave equation \eqref{wave:L} in the spatial
	variable. Finally, the critical H\"older indices of
	$x\mapsto H(t\,,x)$, $x\mapsto W(t\,,x)$ and that of
	the local times of $\bar{X}$ are the same.
\end{theorem}

Blumenthal and Getoor \ycite{BG}
have introduced several ``indices''
that describe various properties of a L\'evy process. 
We recall below their \emph{lower index} $\beta''$:
\begin{equation}
	\beta'' := \liminf_{|\xi|\to\infty} \frac{
	\log\Re\Psi(\xi)}{\log|\xi|} =\sup\left\{
	\alpha\ge 0:\, \lim_{|\xi|\to\infty}
	\frac{\Re\Psi(\xi)}{|\xi|^\alpha} =\infty \right\}.
\end{equation}

\begin{theorem}\label{th:joint:holder}
	If $\beta''>d$, then the stochastic heat equation
	\eqref{heat:L} has
	a random-field solution that is jointly H\"older continuous.
	The critical H\"older index is $\le(\beta''-d)/2$ for the space
	variable and $\le(\beta''-d)/2\beta''$ for the time variable. Furthermore,
	if $X$ is symmetric, then the same assertions hold for the
	solution to the stochastic wave equation \eqref{wave:L}.
\end{theorem}

\begin{remark}
	If $X$ is in the domain of attraction of Brownian motion
	on $\R$, then $\beta''=2$, and the critical temporal and
	spatial index bounds of Theorem \ref{th:joint:holder}
	are respectively $1/2$ and $1/4$. These numbers are
	well known to be the optimal H\"older indices.
	In fact, the relatively simple
	H\"older-index bounds of Theorem \ref{th:joint:holder}
	are frequently sharp; see Example
	\ref{example:stable:Holder}.
	\qed
\end{remark}

It is well known that $0\le \beta''\le 2$
\cite{BG}*{Theorem 5.1}. Thus, the preceding
is inherently a one-dimensional result; this is
in agreement with
the result of Theorem \ref{th:exist}. Alternatively, one
can combine Theorems \ref{th:exist} and \ref{th:joint:holder} 
to construct a proof of the fact
that $\beta''\le 2$. [But of course the original proof 
of Blumenthal and Getoor is simpler.]

Several things need to be made clear here; the first being
the meaning of a ``solution.'' With this aim in mind,
we treat the two equations separately.

We make precise sense of the stochastic
heat equation \eqref{heat:L} in 
much the same manner as Walsh \ycite{Walsh};
see also 
Brze{\'z}niak and van Neerven \ycite{BvN},
Dalang \ycites{Dalang:99,Dalang:01},
and Da Prato \ycite{DaPrato}.
It is well known that it is much harder
to give precise meaning to stochastic hyperbolic equations of
the wave type \eqref{wave:L}, even when $\L$ is the
Laplacian
\cites{DalangLeveque:04a,DalangLeveque:04b,DalangMueller,%
QuerSole,Mueller:97,Dalang:99,Mueller:93}.
Thus, as part of the present work, we introduce
a simple and direct method that makes rigorous sense of \eqref{wave:L}
and other linear SPDEs of this type. We believe
this method to be of some independent interest.

\subsection{The parabolic case}
In order to describe the meaning of \eqref{heat:L}
we need to first introduce some notation. 

Let $\{P_t\}_{t\ge 0}$ denote the semigroup of
the driving L\'evy process $X$; that is,
$(P_tf)(x):=\E f(x+X_t)$ for all bounded Borel-measurable
functions $f:\R^d\to\R$ [say], all $x\in\R^d$, and all $t\ge 0$.
[As usual, $X_0:=0$.] Formally speaking,
$P_t=\exp(t\L )$.

The dual semigroup is denoted by
$P^*$, so that $(P^*_tf)(x) = \E f(x-X_t)$. It is easy to
see that $P$ and $P^*$ are adjoint in $L^2(\R^d)$ in the sense
that
\begin{equation}
	(P_tf\,,g) = (f\,,P_t^*g)
	\quad\text{for all $f,g\in L^2(\R^d)$.}
\end{equation}
Needless to say, $(\cdot\,,\cdot)$ denotes the usual
Hilbertian inner product on $L^2(\R^d)$.

Let $\sS (\R^d)$ denote the class of all rapidly decreasing
test functions on $\R^d$, and recall that:
(i) $\sS (\R^d)
\subset\text{Dom}(\L)\cap C^\infty(\R^d)\subset L^2(\R^d)\cap
C^\infty(\R^d)$; and (ii) for all
$\phi,\psi\in\sS (\R^d)$,
\begin{equation}
	\lim_{h\downarrow 0} \left( \frac{P_{t+h}\phi-P_t\phi}{h}\,,\psi
	\right) = \left( \L  P_t\phi\,,\psi\right).
\end{equation}
Thus, $v(t\,,x):=(P_t\phi)(x)$ solves the Kolmogorov
equation $\partial_t v(t\,,x) = (\L  v)(t\,,x)$,
subject to the initial condition that
$v(0\,,x)=\phi(x)$. This
identifies the Green's function for $\partial_t-\L =0$. Hence,
we can adapt the Green-function method
of Walsh \ycite{Walsh}*{Chapter 3}, without any great
difficulties, to deduce that a \emph{weak solution}
to \eqref{heat:L} is the Gaussian random field
$\{  H (t\,,\phi);\, t\ge 0,\, \phi\in\sS (\R^d)\}$,
where
\begin{equation}\label{heat:L:sol}
	 H (t\,,\phi) := \int_0^t \int_{\R^d}
	(P_{t-s}^*\phi)(y)\, w(dy\,ds).
\end{equation}
This is defined simply as a Wiener integral.

\begin{proposition}\label{pr:heat:exist}
	The Gaussian random field $\{ H (t\,,\phi);\,
	t\ge 0,\, \phi\in\sS (\R^d)\}$ is 
	well defined. Moreover, the process
	$\phi\mapsto  H (t\,,\phi)$ is a.s.\ linear for each
	$t\ge 0$.
\end{proposition}

\begin{proof}
	On one hand,
	the Wiener isometry tells us that
	\begin{equation}
		\E\left( \left|  H (t\,,\phi) \right|^2\right)
		=\int_0^t \|P^*_{t-s}\phi\|_{L^2(\R^d)}^2\, ds.
	\end{equation}
	On the other hand,
	it is known that each $P^*_s$ is a contraction on $L^2(\R^d)$.
	Indeed, it is not hard to check directly that if $\ell_d$ denotes
	the Lebesgue measure on $\R^d$, then the dual L\'evy process
	$-X$ is $\ell_d$-symmetric \cite{FOT}*{pp.\ 27--28}. Therefore,
	the asserted contraction property of $P^*_s$ follows from
	equation (1.4.13) of Fukushima, \={O}shima, and
	Takeda \ycite{FOT}*{p.\ 28}. It follows then that
	$\E ( |  H (t\,,\phi) |^2 )
	\le t\|\phi\|_{L^2(\R^d)}^2$,
	and this is finite for all $t\ge 0$
	and $\phi\in\sS (\R^d)$. This proves that $u$
	is a well-defined Gaussian random field indexed by
	$\R_+\times\sS (\R^d)$. The proof of
	the remaining property follows the argument of Dalang
	\ycite{Dalang:99}*{Section 4} quite
	closely, and is omitted.
\end{proof}

\subsection{The nonrandom hyperbolic case}
It has been known for some time that hyperbolic
SPDEs tend to be harder to study, and even to define
precisely, than their parabolic counterparts. See,
for instance, Dalang and Sanz-Sol\'e
\ycite{DalangSole:07} for the most
recent work on the stochastic wave equation in dimension 3.

In order to define what the stochastic wave equation \eqref{wave:L}
means precisely, we can try to mimic the original Green-function
method of Walsh \ycite{Walsh}. But we quickly run into 
the technical problem of not being able to identify a suitable
Green function (or even a measure) for the corresponding integral
equation. In order to overcome this obstacle,
one could proceed as in Dalang and Sanz-Sol\'e
\ycite{DalangSole:07}, but
generalize the role of their fractional Laplacian. Instead,
we opt for a more direct route that is particularly
well suited for studying the SPDEs of the present type.

In order to understand \eqref{wave:L} better,
we first consider the deterministic integro-differential
equation,
\begin{equation}\label{integraleq}\left|\begin{split}
	&\partial_{tt} u(t\,,x) = (\L  u)(t\,,x) + f(t\,,x),\\
	&u(0\,,x)=\partial_t u(0\,,x)=0,
\end{split}\right.\end{equation}
where $f:\R_+\times\R^d\to\R$ is a ``nice'' function,
and the variables $t$ and $x$ range respectively over
$\R_+$ and $\R^d$. We can study this equation only under
the following symmetry condition:
\begin{equation}\label{assumption:symm}
	\text{The process $X$ is symmetric}.
\end{equation}
Equivalently, we assume that $\Psi$ is \emph{real
and nonnegative}.

Recall that ``\,$\hat{\hskip1mm}$\,'' denotes
the Fourier transform in the
$x$ variable, and apply it informally to \eqref{integraleq}
to deduce that it is equivalent to the following: For all
$t\ge 0$ and $\xi\in\R^d$,
\begin{equation}\left|\begin{split}
	&\partial_{tt} \hat{u}(t\,,\xi) = -\Psi(\xi)
	\hat{u}(t\,,\xi) + \hat{f}(t\,,\xi),\\
	&\hat{u}(0\,,\xi) = \partial_t \hat{u}(0\,,\xi) =0.
\end{split}\right.\end{equation}
This is an inhomogeneous second-order ordinary differential equation
[in $t$] which can be solved explicitly, via Duhamel's principle, 
to produce the following ``formula'':
\begin{equation}
	\hat{u}(t\,,\xi) = \frac{1}{\sqrt{\Psi(\xi)}}
	\int_0^t\sin\left( \sqrt{\Psi(\xi)}\, (t-s) \right)
	\hat{f}(s\,,\xi)\,ds.
\end{equation}
We invert the preceding---informally still---to obtain
\begin{equation}
	u(t\,,x) = \frac{1}{(2\pi)^d}\int_{\R^d}\int_0^t
	\frac{\sin\left(\sqrt{\Psi(\xi)}\, (t-s)\right)}{
	\sqrt{\Psi(\xi)}}\, e^{-i\xi\cdot x}\hat{f}(s\,,\xi)\, ds\, d\xi.
\end{equation}
We can multiply this by a nice function $\phi$, then integrate
$[dx]$ to arrive at 
\begin{equation}\label{wave:def1}
	u(t\,,\phi) := \frac{1}{(2\pi)^d} \int_{
	\R^d}\int_0^t
	\frac{\sin\left(\sqrt{\Psi(\xi)}\, (t-s)\right)}{
	\sqrt{\Psi(\xi)}}\, \overline{\hat\phi(\xi)}\
	\hat{f}(s\,,\xi)\, ds\, d\xi.
\end{equation}
We may think of as this as the
``weak/distributional solution'' to \eqref{integraleq}.

\subsection{The random hyperbolic case}
We follow standard terminology and identify the
white noise $\dot{w}$ with
the iso-Gaussian process $\{w(h)\}_{h\in L^2(\R_+\times\R^d)}$ as follows:
\begin{equation}
	w(h) := \int_{\R^d}\int_0^\infty h(s\,,x)\, w(ds\,dx).
\end{equation}

Next, we define the \emph{Fourier transform} $\hat{w}$ of white noise:
\begin{equation}
	\hat{w}(h) := \frac{w(\hat{h})}{(2\pi)^{d/2}} = \frac{1}{(2\pi)^{d/2}}
	\int_{\R^d}\int_0^\infty 
	\hat{h}(s\,,\xi)\, w(ds\,d\xi),
\end{equation}
all the time remembering that ``$\hat{h}$'' refers
to the Fourier transform of $h$ in its spatial variable. 
Suppose $h(s\,,x)=\phi_1(s)\phi_2(x)$ for $t\ge 0$ and $x\in\R^d$,
where $\phi_1\in L^2(\R_+)$ and
$\phi_2\in L^2(\R^d)$. Then,
it follows from the Wiener isometry and Plancherel's theorem that
\begin{equation}
	\| \hat{w}(h)\|_{L^2(\R_+\times\R^d)}^2=
	\|h\|_{L^2(\R_+\times\R^d)}^2.
\end{equation}
Because $L^2(\R_+)\otimes L^2(\R^d)$ is dense in
$L^2(\R_+\times\R^d)$, this proves that $\hat{w}$ is defined
continuously on all of $L^2(\R_+\times\R^d)$. Moreover,
$\hat{w}$ corresponds to a white noise which is correlated with
$\dot{w}$, as described by the following formula:
\begin{equation}
	\E\left[ w(h_1)\cdot
	\overline{\hat{w}(h_2)}\right]= \frac{1}{(2\pi)^{d/2}}
	\int_0^\infty \int_{\R^d} \overline{\hat{h}_1(s\,,\xi)}\,
	h_2(s\,,\xi)\, d\xi\, ds,
\end{equation}
valid for all $h_1,h_2\in L^2(\R_+\times\R^d)$.

In light of \eqref{wave:def1}, we define the weak solution
to the stochastic wave equation \eqref{wave:L} as the Wiener
integral
\begin{equation}\label{wave:random}
	 W (t\,,\phi) = \frac{1}{(2\pi)^{d/2}}
	\int_{\R^d}\int_0^t \frac{\sin\left(
	\sqrt{\Psi(\xi)}\,(t-s)\right)}{\sqrt{\Psi(\xi)}}\,
	\overline{\hat\phi(\xi)}\, \hat{w}(ds\,d\xi).
\end{equation}
It is possible to verify that this method also works for the
stochastic heat equation \eqref{heat:L}, and that it
produces a equivalent formulation of the Walsh
solution \eqref{heat:L:sol}. However, in the present setting,
this method saves us from having to describe the existence
[and some regularity] of the Green function for
the integral equation \eqref{integraleq}.

\begin{proposition}\label{pr:wave:exist}
	If the symmetry condition
	\eqref{assumption:symm} holds, then the stochastic
	wave equation \eqref{wave:L}
	has a weak solution $ W $ for all $\phi\in\sS (\R^d)$.
	Moreover, $\{ W (t\,,\varphi);\,
	t\ge 0,\, \phi\in\sS (\R^d)\}$ is a well-defined
	Gaussian random field, and 
	$\phi\mapsto  W (t\,,\phi)$ is a.s.\ linear
	for all $t\ge 0$.
\end{proposition}

\begin{proof}
	We apply the Wiener isometry to obtain 
	\begin{equation}
		\E\left(\left|  W (t\,,\phi)\right|^2\right)
		= \frac{1}{(2\pi)^d} \int_{\R^d}\int_0^t 
		\frac{\sin^2\left(
		\sqrt{\Psi(\xi)}\,(t-s)\right)}{\Psi(\xi)}|\hat\phi(\xi)|^2\, dsd\xi.
	\end{equation}
	Because $|\sin\theta/\theta|\le 1$,
	\begin{equation}\begin{split}
		\E\left(\left|  W (t\,,\phi)\right|^2\right)
			&\le\frac{1}{(2\pi)^d} \int_{\R^d}\int_0^t 
			(t-s)^2 |\hat\phi(\xi)|^2\, ds\,d\xi\\
		&=\frac{t^3}{3}\|\phi\|_{L^2(\R^d)}^2,
	\end{split}\end{equation}
	thanks to Plancherel's theorem.
	It follows immediately from
	this that $\{ W (t\,,\phi);\, t\ge
	0,\, \phi\in\sS (\R^d)\}$
	is a well-defined Gaussian random field.
	The remainder of the proposition is standard.
\end{proof}

\section{\bf Existence of functions-valued solutions:
	Proof of Theorem \ref{th:exist}}\label{sec3}
	
Let $u:=\{u(t\,,\phi);\, t\ge 0,\, \phi\in \sS (\R^d)\}$
denote the weak solution to either one of \eqref{heat:L} or
\eqref{wave:L}. Our present goal
is to extend uniquely the Gaussian random field $u$ to 
a Gaussian random field indexed by $\R_+\times M$, 
where $M$ is a maximal subset of $\mathscr{D}(\R^d)$---the
space of all Schwartz distributions on $\R^d$. Such an $M$
exists thanks solely to functional-analytic facts: 
Define, temporarily,
\begin{equation}
	d_t(\phi) := \sqrt{
	\E\left( \left| u(t\,,\phi)\right|^2\right)}\qquad
	\text{for all $\phi\in\sS (\R^d)$
	and $t\ge 0$.}
\end{equation}
Then, the linearity of $u$ in $\phi$ shows that
$(\phi\,,\psi)\mapsto d_t(\phi-\psi)$ defines a metric
for each $t\ge 0$. Let $M_t$ denote the completion
of $\sS (\R^d)$ in $\mathscr{D}(\R^d)$
with respect to the metric
induced by $d_t$; and define $M:=\cap_{t\ge 0}M_t$.
The space
$M$ can be identified with  the
largest possible family of candidate
test functions for weak solutions to either
the stochastic heat equation
\eqref{heat:L} or the stochastic wave
equation \eqref{wave:L}. Standard heuristics from
PDEs then tell us that \eqref{heat:L} and/or \eqref{wave:L}
has random-field solutions if and only if 
$\delta_x\in M$ for all $x\in\R^d$; this can be
interpreted as an equivalent \emph{definition} of random-field solutions.
When $\delta_x\in M$
we may write $u(t\,,x)$ in place of $u(t\,,\delta_x)$.
In order to prove Theorem \ref{th:exist} we will need some
a priori estimates on the weak solutions of both the stochastic
equations \eqref{heat:L} and \eqref{wave:L}. We proceed by
identifying $M$ with generalized Sobolev spaces that arise
in the potential theory of symmetric L\'evy processes.
Now let us begin by studying the parabolic case.

\subsection{The parabolic case}

\begin{proposition}\label{pr:heat:embed}
	Let $ H $ denote the weak solution
	\eqref{heat:L:sol} to the stochastic
	heat equation \eqref{heat:L}.
	Then, for all $\phi\in\sS (\R^d)$,
	$\lambda>0$, and $t\ge 0$,
	\begin{equation}\label{heat:bd}
		\frac{1-e^{-2 t/\lambda}}{2}\sE(\lambda\,;\phi)\le
		\E\left(\left|  H  (t\,,\phi) \right|^2\right)
		\le \frac{e^{2 t/\lambda}}{2}\sE(\lambda\,;\phi),
	\end{equation}
	where
	\begin{equation}\label{def:E_lambda}
		\sE(\lambda\,;\phi) := \frac{1}{(2\pi)^d}
		\int_{\R^d} \frac{|\hat\phi(\xi)|^2}{(1/\lambda)+
		\text{\rm Re}\Psi(\xi)}\, d\xi.
	\end{equation}
\end{proposition}

Next we record the following immediate
but useful corollary; it follows from Proposition
\ref{pr:heat:embed} by simply setting $\lambda:=t$.

\begin{corollary}\label{cor:heat:embed}
	If $ H $ denotes the weak solution to the stochastic
	heat equation \eqref{heat:L},
	then for all $\phi\in\sS (\R^d)$ and $t\ge 0$,
	\begin{equation}
		\tfrac{1}{3}\sE(t\,;\phi)\le
		\E\left(\left|  H  (t\,,\phi) \right|^2\right)
		\le 4\sE(t\,;\phi),
	\end{equation}
\end{corollary}

The preceding upper bound for $\E(|H(t\,,\phi)|^2)$
is closely tied to an energy
inequality for the weakly
asymmetric exclusion process. See Lemma 3.1
of \ocite{BertiniGiacomin};
they ascribe that lemma to H.-T. Yau.

The key step of the proof of Proposition
\ref{pr:heat:embed} is an elementary real-variable
result which we prove next.
\begin{lemma}\label{lem:Monotone}
	If $g:\R_+\to\R_+$ is Borel measurable and
	nonincreasing, then
	for all $t,\lambda>0$,
	\begin{equation}
		\left(1-e^{-2t/\lambda}\right)\,
		\int_0^\infty e^{-2s/\lambda}g(s)\, ds
		\le \int_0^t g(s)\, ds\\
		\le e^{2t/\lambda}  \int_0^\infty e^{-2s/\lambda}g(s)\, ds.
	\end{equation} 
	Monotonicity is not needed for the upper bound on
	$\int_0^tg(s)\,ds$.
\end{lemma}

\begin{proof}[Proof of Lemma \ref{lem:Monotone}]
	The upper bound on $\int_0^t g(s)\,ds$ follows simply because
	$e^{2(t-s)/\lambda}\geq 1$ whenever $t\ge s$. 
	In order to derive the lower bound we write
	\begin{equation}\begin{split}
		\int_0^\infty e^{-2 s/\lambda} g(s)\, ds 
			&=\sum_{n=0}^\infty \int_{nt}^{(n+1)t}
			e^{-2s/\lambda} g(s)\, ds\\
		&\le \sum_{n=0}^\infty e^{-2nt/\lambda} \int_0^t g(s+nt)\, ds.
	\end{split}\end{equation}
	Because $g$ is nonincreasing we can write
	$g(s+nt)\le g(s)$ to conclude the proof.
\end{proof}

\begin{proof}[Proof of Proposition \ref{pr:heat:embed}]
	We know from \eqref{heat:L:sol} and the Wiener isometry
	that
	\begin{equation}\label{eq:3.7}
		\E\left(\left|  H (t\,,\phi)\right|^2\right)
		=\int_0^t \| P^*_s \phi\|_{L^2(\R^d)}^2\, ds.
	\end{equation}
	Since the Fourier multiplier of $P^*_s$
	is $\exp(-s\Psi(-\xi))$ at $\xi\in\R^d$,
	we can apply the Plancherel theorem 
	and deduce the following formula:
	\begin{equation}\label{eq:3.8}
		\| P^*_s \phi\|_{L^2(\R^d)}^2 =
		\frac{1}{(2\pi)^d} \int_{\R^d} e^{-2s
		\text{\rm Re}\Psi(\xi)}|\hat\phi(\xi)|^2\, d\xi.
	\end{equation}
	Because $\text{\rm Re}\Psi(\xi)\ge 0$, Lemma \ref{lem:Monotone}
	readily proves the proposition.
\end{proof}

Equation \eqref{def:E_lambda} can be used to define 
$\sE(\lambda\,;\phi)$ for all Schwartz distributions
$\phi$, and not only those in $\sS (\R^d)$. Moreover,
it is possible to verify directly that $\phi\mapsto
\sE(\lambda\,;\phi)^{1/2}$ defines a norm on $\sS (\R^d)$.
But in all but uninteresting cases, $\sS (\R^d)$
is \emph{not} complete in this norm. Let
$L^2_\L(\R^d)$ denote the completion of $\sS (\R^d)$
in the norm $\sE(\lambda\,;\bullet)^{1/2}$. Thus the Hilbert space 
$L^2_\L(\R^d)$ can be identified with $M$. The following is a result
about the potential theory of symmetric L\'evy processes, but
we present a  self-contained proof that does not depend on
that deep theory.

\begin{lemma}\label{lem:Sobolev}
	The space $L^2_\L(\R^d)$ does not depend on
	the value of $\lambda$. Moreover, $L^2_\L(\R^d)$
	is a Hilbert space in norm
	$\sE(\lambda\,;\bullet)^{1/2}$ for each fixed $\lambda>0$. Finally,
	the quasi-isometry
	\eqref{heat:bd} is valid for all $t\ge 0$, $\lambda>0$,
	and $\phi\in L^2_\L(\R^d)$.
\end{lemma}

\begin{proof}
	We write, temporarily, $L^2_{\L,\lambda}(\R^d)$
	for $L^2_\L(\R^d)$, and seek to prove that it is
	independent of the choice of $\lambda$.
	
	Define for all distributions $\phi$ and $\psi$,
	\begin{equation}
		\sE(\lambda\,;\phi\,,\psi) :=
		\frac{1}{2(2\pi)^d}\left[
		\int_{\R^d} \frac{
		\hat\phi(\xi) \ \overline{\hat\psi(\xi)}}{
		(1/\lambda)+\text{\rm Re}\Psi(\xi)}\, d\xi
		+
		\int_{\R^d} \frac{
		\hat\psi(\xi) \ \overline{\hat\phi(\xi)}}{
		(1/\lambda)+\text{\rm Re}\Psi(\xi)}\, d\xi
		\right].
	\end{equation}
	For each $\lambda>0$ fixed,
	$(\phi\,,\psi)\mapsto\sE(\lambda\,;\phi\,,\psi)$
	is a pre-Hilbertian inner product on
	$\sS (\R^d)$, and $\sE(\lambda\,;\phi)
	=\sE(\lambda\,;\phi\,,\phi)$.
	
	Thanks to Proposition \ref{pr:heat:embed}, for
	all $\alpha>0$ there exists a finite and positive
	constant $c=c_{\alpha,\lambda}$ such that
	$c^{-1}\sE(\alpha\,;\phi)
	\le \sE(\lambda\,;\phi) \le c\sE(\alpha\,;\phi)$
	for all $\phi\in\sS (\R^d)$. This proves that
	$L^2_{\L ,\lambda}(\R^d)=L^2_{\L ,\alpha}(\R^d)$,
	whence follows the independence of $L^2_\L(\R^d)$
	from the value of $\lambda$. The remainder of the lemma
	is elementary.
\end{proof}

The space $L^2_\L(\R^d)$ is a generalized
Sobolev space, and contains many classical spaces of
Bessel potentials, as the following example shows.

\begin{example}
	Suppose $\L =-(-\Delta)^{s/2}$ for 
	$s\in(0\,,2]$. Then, $\L $ is the generator
	of an isotropic stable-$s$ L\'evy process, and
	$L^2_\L(\R^d)$ is the space $H_{-s/2}(\R^d)$
	of Bessel potentials. For a nice pedagogic treatment
	see the book of Folland \ycite{Folland}*{Chapter 6}.\qed
\end{example}

\subsection{The hyperbolic case}
The main result of this section is the following
quasi-isometry; it is the wave-equation
analogue of Proposition \ref{pr:heat:embed}.

\begin{proposition}\label{pr:wave:embed}
	Suppose the symmetry condition \eqref{assumption:symm}
	holds, and let $ W :=\{ W 
	(t\,,\phi);\, t\ge 0,\, \phi\in\sS (\R^d)\}$
	denote the weak solution to the stochastic wave
	equation \eqref{wave:L}. Then,
	\begin{equation}
		\tfrac14 t \sE\left( t^2;\phi \right)\le
		\E\left(\left|  W  (t\,,\phi)\right|^2\right) 
		\le 2t\sE \left( t^2;\phi \right).
	\end{equation}
	for all $t\ge 0$ and $\phi\in\sS (\R^d)$.
	Moreover,
	we can extend $ W $ by density so that the preceding
	display continues to remain valid 
	when $t\ge 0$ and $\phi\in L^2_\L(\R^d)$.
\end{proposition}

\begin{proof}
	Although Lemma \ref{lem:Monotone} is not applicable,
	we can proceed in a similar manner as we did when
	we proved the earlier quasi-isometry
	result for the heat equation (Proposition \ref{pr:heat:embed}).
	Namely, we begin by observing that
	\begin{equation}\label{formula:wave}
		\E\left(\left|  W (t\,,\phi) \right|^2\right) =
		\frac{1}{(2\pi)^d} \int_0^t \int_{\R^d}
			\frac{\sin^2\left(\sqrt{\Psi(\xi)}\, s\right)}{
			\Psi(\xi)}
			\ | \hat{\phi}(\xi) |^2\, d\xi\, ds.
	\end{equation}
	See \eqref{wave:random}.
	If $\theta>0$ then $\sin\theta$ is at most the minimum of
	one and $\theta$. This leads to
	the bounds
	\begin{equation}\begin{split}
		\E\left(\left|  W (t\,,\phi) \right|^2\right) &\le
			\frac{1}{(2\pi)^d} \int_0^t \int_{\R^d}
			\left(s^2\wedge \frac{1}{\Psi(\xi)}
			\right)
			\ | \hat{\phi}(\xi) |^2\, d\xi\, ds\\
		&\le \frac{t}{(2\pi)^d} \int_{\R^d}
			\left(t^2\wedge \frac{1}{\Psi(\xi)}
			\right)
			\ | \hat{\phi}(\xi) |^2 \, d\xi.
	\end{split}\end{equation}
	The upper
	bound follows from this and the elementary inequality
	$t^2\wedge z^{-1}\le 2/(t^{-2}+z)$, valid for all $z\ge 0$.
	
	In order to derive the [slightly] harder lower bound
	we first rewrite \eqref{formula:wave} as follows:
	\begin{equation}\label{wave:Eq0}
		\E\left(\left|  W (t\,,\phi) \right|^2\right) =
		\frac{t}{2(2\pi)^d} \int_{\R^d}
		\left( 1 - \frac{\sin\left( 2\sqrt{\Psi(\xi)}\, t\right)}{
		2\sqrt{\Psi(\xi)}\, t}\right)
		\frac{\ | \hat{\phi}(\xi) |^2}{\Psi(\xi)}\, d\xi.
	\end{equation}
	We shall analyze the integral by splitting it according
	to whether or not $\Psi\le 1/t^2$.
	
	Taylor's expansion [with remainder] reveals that
	if $\theta$ is nonnegative, then $\sin \theta$
	is at most $\theta - (\theta^3/6) + (\theta^5/120)$. This
	and a little algebra together show that
	\begin{equation}
		1- \frac{\sin\theta}{\theta} \ge \frac{2\theta^2}{15}
		\qquad\text{if $0\le\theta\le 2$}.
	\end{equation}
	Consequently,
	\begin{equation}\label{wave:Eq1}\begin{split}
		\int_{\{\Psi\le 1/t^2\}}\left(
			1- \frac{\sin\left(2\sqrt{\Psi(\xi)}\, t\right)}{
			2\sqrt{\Psi(\xi)}\, t}\right) \frac{|\hat\phi(\xi)|^2}{
			\Psi(\xi)}\, d\xi
			&\ge \frac{8t^2}{15}\int_{\{\Psi\le 1/t^2\}}
			|\hat\phi(\xi)|^2\, d\xi\\
		&\ge\frac{1}{2}\int_{\{\Psi\le 1/t^2\}}
			\left( t^2 \wedge \frac{1}{\Psi(\xi)}
			\right) |\hat\phi(\xi)|^2 \, d\xi.
	\end{split}\end{equation}
	
	For the remaining integral we use the elementary bound
	$1-(\sin\theta/\theta)\ge 1/2$, valid for all $\theta> 2$.
	This leads to the following inequalities:
	\begin{equation}\label{wave:Eq2}\begin{split}
		\int_{\{\Psi> 1/t^2\}}\left(
			1- \frac{\sin\left(2\sqrt{\Psi(\xi)}\, t\right)}{
			2\sqrt{\Psi(\xi)}\, t}\right) \frac{|\hat\phi(\xi)|^2}{
			\Psi(\xi)}\, d\xi
			&\ge \frac{1}{2}\int_{\{\Psi> 1/t^2\}}
			\frac{|\hat\phi(\xi)|^2}{\Psi(\xi)}\, d\xi\\
		&= \frac{1}{2}\int_{\{\Psi> 1/t^2\}}
			\left( t^2\wedge \frac{1}{\Psi(\xi)}\right)
			|\hat\phi(\xi)|^2 \, d\xi.
	\end{split}\end{equation}
	The proof concludes from summing up 
	equations \eqref{wave:Eq1} and \eqref{wave:Eq2},
	and then plugging the end result into \eqref{wave:Eq0}.
\end{proof}
We now give a proof of Theorem  \ref{th:exist}.
\begin{proof}
	Let us begin with the proof in the case of the stochastic
	heat equation \eqref{heat:L}.
	Proposition \ref{pr:heat:embed} is a quasi-isometry
	of the maximal space $M$ of test functions
	for weak solutions of \eqref{heat:L} into
	$L^2_\L(\R^d)$. Thus, $M$ can be
	identified with the Hilbert space $L^2_\L(\R^d)$,
	and hence
	\eqref{heat:L} has random-field
	solutions if and only if $\delta_x\in L^2_\L(\R^d)$
	for all $x\in\R^d$.  Thanks to Lemma \ref{lem:Sobolev},
	the stochastic heat equation
	\eqref{heat:L} has random-field
	solutions if and only if
	\begin{equation}\label{cond:hawkes}
		\int_{\R^d} \frac{d\xi}{\vartheta+\text{\rm Re}\Psi(\xi)}< \infty
		\qquad\text{ for some, 
		and hence all, $\vartheta>0$.}
	\end{equation}
	The first part of the proof is concluded since
	condition \eqref{cond:hawkes} is known to be necessary 
	as well as sufficient for $\bar{X}$ to have local times
	\cite{Hawkes}*{Theorem 1}. The remaining portion of the
	proof follows from the preceding in much the
	same way as the first portion was deduced from
	Proposition \ref{pr:heat:embed}.
\end{proof}

\section{\bf Spatial continuity: Proof of Theorem \ref{th:cont}}\label{sec4}
\begin{proof}We work with the stochastic heat equation
	\eqref{heat:L} first. Without loss of generality, we may---and will---assume
	that \eqref{heat:L} has a random-field solution $ H (t\,,x)$,
	and $\bar{X}$ has local times. Else, Theorem \ref{th:exist} finishes the proof.
	
	Let $\phi:=\delta_x-\delta_y$, 
	and note that $|\hat\phi(\xi)|^2=2(1-\cos(\xi(x-y)))$ is a function
	of $x-y$. Because
	$ H (t\,,\phi)= H (t\,,x)- H (t\,,y)$,
	equations \eqref{eq:3.7} and \eqref{eq:3.8} imply that
	$z\mapsto  H (t\,,z)$ is a centered Gaussian process with
	stationary increments for each fixed $t\ge 0$.
	
	Consider the function
	\begin{equation}\label{eq:Omega}
		 h (r) := \frac{1}{2\pi}\int_{-\infty}^\infty
		\frac{1-\cos(r\xi)}{1+\Re \Psi(\xi)}\, d\xi,\qquad
		\text{defined for all $r\ge 0$}.
	\end{equation}
	
	According to Lemma \ref{lem:Sobolev},
	Proposition \ref{pr:heat:embed} holds for all
	Schwartz distributions $\phi\in L^2_\L(\R^d)$.
	The existence of random-field solutions is equivalent
	to the condition that $\delta_x\in L^2_\L(\R^d)$
	for all $x\in\R$. We apply Proposition \ref{pr:heat:embed}
	to $\phi:=\delta_x-\delta_y$, with $\lambda:=1$ [say],
	and find that
	\begin{equation}\label{eq:hEh}
		\left( 1- e^{-2t}\right)  h (|x-y|)\le
		\E\left(\left|  H (t\,,x) -  H (t\,,y) \right|^2\right)
		\le e^{2t}  h (|x-y|).
	\end{equation}
	Define $\bar{ h }$ to be the Hardy--Littlewood nondecreasing
	rearrangement of $ h $. That is,
	\begin{equation}\label{monotone:rearrangement}
		\bar{ h }(r) := \inf\{ y\ge 0:\, g  (y) > r\}
		\quad\text{where}\quad
		g (y) := \text{meas}\left\{ r\ge 0:\,  h (r) \le y\right\}.
	\end{equation}
	Then according to the proof of Corollary 6.4.4 of
	Marcus and Rosen
	\ycite{MarcusRosen}*{p.\ 274}, the stationary-increments
	Gaussian process
	$x\mapsto  H (t\,,x)$ has a continuous modification iff
	\begin{equation}\label{cond:Barlow}
		\int_{0^+} \frac{\bar{ h }(r)}{r
		| \log r |^{1/2}}\, dr<\infty.
	\end{equation}
	Next we claim that \eqref{cond:Barlow} is equivalent to 
	the continuity of the local times of the symmetrized L\'evy
	process $\bar{X}$. We recall that the characteristic exponent
	of $\bar{X}$ is $2\Re\Psi$, and hence by the L\'evy--Khintchine
	formula it can be written as
	\begin{equation}
		2\Re\Psi(\xi) = \sigma^2 \xi^2 + \int_{-\infty}^\infty
		( 1-\cos(\xi x))\, \nu(dx),
	\end{equation}
	where $\nu$ is a $\sigma$-finite Borel measure on $\R$
	with $\int_{-\infty}^\infty (1\wedge x^2)\, \nu(dx)<\infty$.
	See, for example Bertoin
	\ycite{Bertoin}*{Theorem 1, p.\ 13}.
	
	Suppose, first, that 
	\begin{equation}\label{eq:LT:E}
		\text{either $\sigma^2>0$ \ or $\int_{-\infty}^\infty
		(1\wedge |x|)\,\nu(dx)=\infty$.}
	\end{equation}
	Then, \eqref{cond:Barlow} is also necessary and sufficient
	for the [joint] continuity of the local times
	of the symmetrized process $\bar{X}$;
	confer with Barlow \ycite{Barlow}*{Theorems B and 1}.
	On the other hand, if \eqref{eq:LT:E} fails to hold then
	$\bar{X}$ is a compound Poisson process.
	Because $\bar{X}$ is also a symmetric process, the L\'evy--Khintchine
	formula tells us that it has zero drift. That is,
	$\bar{X}_t=\sum_{j=1}^{\Pi(t)} Z_i$, where $\{Z_i\}_{i=1}^\infty$
	are i.i.d.\ and symmetric, and $\Pi$ is an independent Poisson process.
	It follows immediately from this that the range of $\bar{X}$
	is a.s.\ countable in that case. This proves that the occupation measure for
	$\bar{X}$ is a.s.\ singular with respect to Lebesgue measure, and therefore
	$\bar{X}$ cannot possess local times. But
	that contradicts the original
	assumption that $\bar{X}$ has local times. Consequently, \eqref{cond:Barlow}
	implies \eqref{eq:LT:E},
	and is equivalent to the spatial continuity of [a modification of] local times
	of $\bar{X}$. This proves the theorem in the parabolic case.
	
	Let us assume further the symmetry condition
	\eqref{assumption:symm}. Because
	$\sE(t^2;\phi)/\sE(1\,;\phi)$ is bounded
	above and below by positive finite constants that depend only on
	$t>0$ (Proposition \ref{pr:heat:embed}), we can apply the
	very same argument to the stochastic wave equation
	\eqref{wave:L}, but use Proposition \ref{pr:wave:embed}
	in place of Proposition \ref{pr:heat:embed}.
	This completes our proof of Theorem \ref{th:cont}.
	\end{proof}

\section{\bf An aside on temporal continuity}\label{sec5}
We spend a few pages discussing matters of
temporal continuity---especially temporal
H\"older continuity---of weak solutions of the stochastic heat
equation \eqref{heat:L}, as well as the stochastic wave
equation \eqref{wave:L}.  

\begin{definition}\label{def:g}
	We call a function $g(s)$ a {\it gauge function}
	if the following are satisfied:
	\begin{enumerate}
		\item $g(s)$ is an increasing function;
		\item $g(s)$ is a slowly varying function at infinity;
		\item $g(s)$ satisfies the integrability condition,
			\begin{equation}
				\int_{0^+} \frac{ds}{s\log(1/s)g(1/s)}<\infty.
			\end{equation}
	\end{enumerate}
\end{definition}

Next, we quote a useful property of slow varying functions
\cite{BGT}*{p.\ 27}.

\begin{proposition}\label{prop:slow}
	If $g$ is a slowly varying function and $\alpha >1$,
	then the integral $\int_x^\infty t^{-\alpha} g(t)\,dt$ converges 
	for every $x>0$, and
	\begin{equation}
		\int_x^\infty \frac{g(t)}{t^\alpha}\,dt
		\sim \frac{g(x)}{(\alpha-1)x^{\alpha-1}}
		\qquad\text{as $x\to\infty$}.
	\end{equation}
\end{proposition}

We can now state the main theorem of this section. It gives a criteria for the temporal continuity of the weak solutions of our stochastic equations.

\begin{theorem}\label{th:cont:t}
	Let $ H $ denote the weak solution to the stochastic 
	heat equation \eqref{heat:L}.  Let $g$ be a gauge
	function in the sense of Definition \ref{def:g}. Choose and fix
	$\phi\in L^2_\L(\R^d)$. Then, $t\mapsto  H (t\,,\phi)$
	has a continuous modification if  the following is satisfied
	\begin{equation}\label{eq:cont:t}
		\int_{\R^d}\frac{\log(1+\Re \Psi(\xi))
		g(1+ |\Psi(\xi)|)}{1+\Re \Psi(\xi)}
		|\hat{\phi}(\xi)|^2
		\,d\xi<\infty.
	\end{equation}
	Moreover, the critical H\"older exponent of $t\mapsto 
	 H (t\,,\phi)$ is precisely
	\begin{equation}\label{eq:Hcont:t}
		\underline{\text{\em ind}}\, \sE(\bullet\,;\phi) :=
		\liminf_{\epsilon\downarrow 0}\frac{\log\sE
		(\epsilon\,;\phi)}{\log\epsilon}.
	\end{equation}
	Consequently, $t\mapsto  H (t\,,\phi)$
	has a H\"older-continuous modification (a.s.) iff
	$\underline{\text{\rm ind}}\,\sE(\bullet\,;\phi)>0$.
	
	If the symmetry condition \eqref{assumption:symm} holds,
	then \eqref{eq:cont:t} guarantees the existence of
	a continuous modification of $t\mapsto  W (t\,,\phi)$,
	where $ W $
	denotes the weak solution to the stochastic wave equation.
	Furthermore, \eqref{eq:Hcont:t} implies the temporal 
	H\"older continuity of $t\mapsto  W (t\,,\phi)$
	of any order
	$<\frac12\, \underline{\text{\em ind}}\, \sE(\bullet\,;\phi)$.
\end{theorem}

We now give two examples.  The first one is about the temporal H\"older exponent while the second one provides a family of random-field solutions which are almost surely unbounded in every open space-time set.
\begin{example}\label{example:stable:Holder}
	Suppose $d=1$
	and $\L =-(-\Delta)^{\alpha/2}$ for some $\alpha\in(1\,,2]$.
	Choose and fix $x\in\R$.
	According to Theorem \ref{th:exist}, $\delta_x\in L^2_\L (\R)$,
	so we can apply Theorem \ref{th:cont:t} with $\phi:=\delta_x$.
	In this case,
	\begin{equation}\begin{split}
		\sE(\epsilon\,;\delta_x) &= \frac{1}{2\pi}
			\int_{-\infty}^\infty \frac{d\xi}{(1/\epsilon)+
			|\xi|^\alpha}\\
		&=\text{const}\cdot\epsilon^{1-(1/\alpha)}.
	\end{split}\end{equation}
	In particular, $\underline{\text{\rm ind}}\,
	\sE(\bullet\,;\delta_x)=1-(1/\alpha)$, whence it follows
	that the critical H\"older exponent of $t\mapsto
	 H (t\,,x)$ is precisely $\frac12 - (2\alpha)^{-1}$.
	When $\L =\Delta$, we have $\alpha=2$ and the critical
	temporal exponent is $1/4$, which agrees with a well-known
	folklore theorem. For example, Corollary 3.4 of Walsh
	\ycite{Walsh}*{pp.\ 318--320}
	and its proof contain this statement for the closely-related
	stochastic cable equation. In particular, see the last two lines
	on page 319 of Walsh's lectures (\emph{loc.\ cit.}).\qed
\end{example}

\begin{example}\label{example:sharp}
	Choose and fix $\alpha\in\R$.
	According to the L\'evy--Khintchine formula
	\cite{Bertoin}*{Theorem 1, p.\ 13} we can find 
	a symmetric L\'evy process whose characteristic exponents satisfies
	$\Psi(\xi) \sim \xi( \log \xi)^\alpha$
	as $\xi\to\infty$.
	Throughout, we assume that $\alpha>1$. This ensures that 
	condition \eqref{cond:hawkes} is in place; i.e.,
	$(1+ \Psi)^{-1}\in L^1(\R)$. Equivalently, that
	both SPDEs \eqref{heat:L} and \eqref{wave:L} have a
	random-field solution.
	
	A few lines of computations show that for all $x\in\R$,
	\begin{equation}\label{eq:example:cont:t}
		\sE(\epsilon\,;\delta_x) \asymp \left| \log(1/\epsilon)
		\right|^{-\alpha+1},
	\end{equation}
	where $a(\epsilon)\asymp b(\epsilon)$ means that $a(\epsilon)/
	b(\epsilon)$ is bounded
	and  below by absolute constants, uniformly for all $\epsilon>0$
	sufficiently small and $x\in\R$. Next, consider
	\begin{equation}
		d(s\,,t) := \sqrt{\E\left(\left|  H (s\,,x)
		- H (t\,,x) \right|^2\right)}.
	\end{equation}
	According to Proposition \ref{pr:mod-in-t:linear} below,
	\eqref{eq:example:cont:t} implies that
	\begin{equation}
		d(s\,,t) \asymp \left| \log \left( \frac{1}{|s-t|}
		\right) \right|^{(1-\alpha)/2},
	\end{equation}
	uniformly for all $s$ and $t$ in a fixed compact subset 
	$[0\,,T]$ of $\R_+$, say. Let $N$ denote the \emph{metric entropy}
	of $[0\,,T]$ in the [pseudo-] metric $d$
	\cite{Dudley}. That is,
	for all $\epsilon>0$, we
	define $N(\epsilon)$ to be the minimum number
	of $d$-balls of radius $\epsilon$
	needed to cover $[0\,,T]$. Then, it is easy
	to deduce from the previous display that
	$\log N(\epsilon) \asymp\epsilon^{-2/(\alpha-1)}$, and hence
	\begin{equation}
		\lim_{\epsilon\downarrow 0} \epsilon\sqrt{\log N(\epsilon)}
		=0\quad\text{if and only if}\quad \alpha>2.
	\end{equation}
	Therefore, if $1<\alpha\le 2$ then Sudakov minorization
	\cite{MarcusRosen}*{p.\ 250} tells us that the process $t\mapsto
	H (t\,,x)$---and also $t\mapsto W (t\,,x)$---does \emph{not}
	have any continuous modifications, almost surely.
	On the other hand, if $\alpha>2$, then 
	the integrability condition \eqref{eq:cont:t}
	holds manifestly, and hence $t\mapsto  H (t\,,x)$
	and $t\mapsto  W (t\,,x)$ both have continuous
	modifications. Thus, the sufficiency condition 
	of Theorem \ref{th:cont:t} is also necessary
	for the present example. In addition, when $\alpha\le 2$,
	the random-field solutions $H$ and $W$ are both unbounded
	a.s.\ in every open set. This assertion follows from general
	facts about Gaussian processes; see, for example equation
	(6.32) of Marcus and Rosen
	\ycite{MarcusRosen}*{p.\ 250}. We skip the details.
	\qed
\end{example}

\subsection{Estimate in the parabolic case}

\begin{proposition}\label{pr:mod-in-t:linear}
	For every $t,\epsilon\ge 0$ and
	$\phi\in L^2_\L (\R^d)$,
	\begin{equation}
		\tfrac12\,  \sE ( \epsilon\,;\phi ) \le
		\E\left( \left|  H (t+\epsilon\,,\phi)
		-  H (t\,,\phi)\right|^2 \right)\le \sE ( \epsilon\,;\phi )
		+ e^{2t} \mathscr{F}(\epsilon\,;\phi),
	\end{equation}
	where
	\begin{equation}\label{eq:F}
		\mathscr{F}(\epsilon\,;\phi) :=
		\frac{1}{(2\pi)^d} \int_{\R^d}\left(1\wedge
		\epsilon^2|\Psi(\xi)|^2\right)
		\frac{|\hat\phi(\xi)|^2}{1+\Re\Psi(\xi)}\,d\xi.
	\end{equation}
\end{proposition}

\begin{proof}
	Because of density it suffices
	to prove that the proposition holds for all functions
	$\phi\in \sS(\R^d)$ of rapid decrease.
	By \eqref{heat:L:sol} and Wiener's isometry,
	\begin{equation}\begin{split}
		\E\left(\left|  H (t+\epsilon\,,\phi) -  H (t\,,\phi)
			\right|^2 \right)
			&= \int_0^t \int_{\R^d} \left|
			\left( P_{t-s+\epsilon}^* \phi \right) (y) -
			\left( P_{t-s}^* \phi\right) (y)
			\right|^2\, dy \, ds\\
		&\hskip1.14in + \int_t^{t+\epsilon}
			\int_{\R^d} \left| \left(
			P_{t-s+\epsilon}^* \phi\right)(y)\right|^2\, dy\, ds.
	\end{split}\end{equation}
	We apply Plancherel's theorem to find that
	\begin{equation}\begin{split}
		\E\left(\left|  H (t+\epsilon\,,\phi) -  H (t\,,\phi)
			\right|^2 \right)
			&= \frac{1}{(2\pi)^d} \int_0^t \int_{\R^d} \left|
			e^{-(s+\epsilon)\Psi(-\xi)} - e^{-s\Psi(-\xi)}
			\right|^2\,  | \hat\phi(\xi) |^2 \, d\xi \, ds\\
		&\hskip0.85in + \frac{1}{(2\pi)^d} \int_0^\epsilon
			\int_{\R^d} e^{-2s\Re\Psi(\xi)}
			\,  | \hat\phi(\xi) |^2 \, dy\, ds.
	\end{split}\end{equation}
	Thus, we can write
	\begin{equation}
		\E\left(\left|  H (t+\epsilon\,,\phi) -  H (t\,,\phi)
		\right|^2 \right) := \frac{T_1 +T_2}{(2\pi)^d},
	\end{equation}
	where
	\begin{align}
		T_1 &:= \int_0^t \int_{\R^d} e^{-2s\Re\Psi(\xi)}
			\left| 1-e^{-\epsilon\Psi(\xi)}\right|^2
			|\hat\phi(\xi)  |^2\, d\xi\, ds,\\
		\intertext{and}
		T_2 &:= \int_0^\epsilon\int_{\R^d} e^{-2s\Re\Psi(\xi)}
			| \hat\phi(\xi) |^2\, d\xi\, ds.
	\end{align}
	
	First we estimate $T_2$, viz.,
	\begin{equation}
		\int_0^\epsilon e^{-2s\Re\Psi(\xi)}\, ds =
		\epsilon\, \frac{1-e^{-2\epsilon\Re\Psi(\xi)}}{2\epsilon
		\Re\Psi(\xi)}.
	\end{equation}
	Because
	\begin{equation}\label{var:theta}
		\frac12\, \frac{1}{1+\theta} \le
		\frac{1-e^{-\theta}}{\theta} \le \frac{2}{1+\theta}
		\quad\text{for all}\  \theta>0,
	\end{equation}
	it follows that
	$\tfrac12(2\pi)^d\,  \sE(\epsilon\,;\phi) \le
	T_2 \le (2\pi)^d\sE(\epsilon\,;\phi).$
	Since $T_1\ge 0$ we obtain the first inequality of
	the proposition.
	
	According to Lemma \ref{lem:Monotone} with $\lambda=1$,
	\begin{equation}
		T_1 \le \frac{e^{2t}}{2}
		\int_{\R^d} \left| 1 - e^{-\epsilon\Psi(\xi)}\right|^2
		\, \frac{| \hat\phi(\xi) |^2}{1+\Re\Psi(\xi)}
		\, d\xi.
	\end{equation}
	Because $| 1 - e^{-\epsilon\Psi(\xi)} |^2 
	\le 1\wedge \epsilon^2 |\Psi(\xi) |^2$, it follows that 
	$T_1\le (2\pi)^d \exp(2t) \F(\epsilon\,;\phi)$, and 
	hence the proof is completed.
\end{proof}

\subsection{Estimate in the hyperbolic case}

\begin{proposition}\label{pr:wave:t-est}
	Assume the symmetry condition \eqref{assumption:symm},
	and let $ W $ denote the weak solution to the stochastic
	wave equation \eqref{wave:L}. Then,
	for all $t\ge 0$, $\epsilon>0$, and $\phi\in
	L^2_\L (\R^d)$,
	\begin{equation}
		\E\left(\left|  W (t+\epsilon\,,\phi) -  W (t\,,\phi)
		\right|^2\right) \le (8t+6\epsilon)\,\sE ( \epsilon^2;
		\phi ).
	\end{equation}
\end{proposition}

\begin{proof}
	By density, it suffices to prove the proposition
	for all functions
	$\phi\in\sS(\R^d)$ of rapid decrease. 
	Henceforth, we choose
	and fix such a function $\phi$.
	
	In accord with \eqref{wave:random} we write
	\begin{equation}
		\E\left(\left|  W (t+\epsilon\,,\phi) -  W (t\,,\phi)
		\right|^2\right) := T_1+T_2,
	\end{equation}
	where
	\begin{equation}
		T_1
		:=\frac{1}{(2\pi)^d}\int_0^t\int_{\R^d}
		\frac{\left|\sin\left(\sqrt{\Psi(\xi)}\,
		(r+\epsilon) \right) - \sin\left(
		\sqrt{\Psi(\xi)}\, r\right) \right|^2}{
		\Psi(\xi)} |\hat\phi(\xi)|^2\, d\xi\, dr,
	\end{equation}
	and
	\begin{equation}
		T_2 :=\frac{1}{(2\pi)^d} \int_t^{t+\epsilon}
		\int_{\R^d} \frac{\sin^2\left(
		\sqrt{\Psi(\xi)}\, r\right)}{\Psi(\xi)}\,
		|\hat\phi(\xi)|^2\, d\xi\, dr.
	\end{equation}
	
	We estimate $T_2$ first:
	The argument that led to \eqref{wave:Eq2} also leads to
	the following inequality:
	\begin{equation}
		\int_{\R^d} \frac{\sin^2\left(
		\sqrt{\Psi(\xi)}\, r\right)}{\Psi(\xi)}\,
		|\hat\phi(\xi)|^2\, d\xi
		\le \int_{\R^d} \left( 
		r^2 \wedge \frac{1}{\Psi(\xi)}\right)|\hat\phi(\xi)|^2
		\, d\xi.
	\end{equation}
	Because
	$\int_\epsilon^{2\epsilon} \left( r^2\wedge a \right) \, dr
	\le 3\epsilon\left( \epsilon^2\wedge a\right)$ for all $a,\epsilon>0$,
	\begin{equation}\label{estimate:T2}\begin{split}
		T_2 &\le \frac{3\epsilon}{(2\pi)^d} \int_{\R^d}
			\left( \epsilon^2\wedge \frac{1}{\Psi(\xi)}\right)
			|\hat\phi(\xi)|^2\, d\xi\\
		&\le 6\epsilon\,\sE ( \epsilon^2;\phi ).
	\end{split}\end{equation}
	
	The estimate for $T_1$ is even simpler to
	derive: Because $|\sin\alpha-\sin\beta |^2$ is bounded above by
	the minimum of $4$ and $2[1-\cos(\beta-\alpha)]$,
	\begin{equation}\begin{split}
		&\int_{\R^d}
			\frac{\left|\sin\left(\sqrt{\Psi(\xi)}\,
			(r+\epsilon) \right) - \sin\left(
			\sqrt{\Psi(\xi)}\, r\right) \right|^2}{
			\Psi(\xi)} |\hat\phi(\xi)|^2\, d\xi\\
		&\hskip2.7in \le \int_{\R^d} \frac{
			1-\cos\left( \sqrt{\Psi(\xi)}\, \epsilon\right)}{\Psi(\xi)}
			|\hat\phi(\xi)|^2\, d\xi.
	\end{split}\end{equation}
	This and the elementary inequality
	$1-\cos x\le x^2/2$ together yield the bound
	\begin{equation}
		\frac{1-\cos\left(\sqrt{\Psi(\xi)}\, \epsilon\right)}{
		\Psi(\xi)}\le 2(\epsilon^2\wedge \frac{1}{\Psi(\xi)}).
	\end{equation}
	Consequently,
	\begin{equation}
		\frac{1}{(2\pi)^d}\int_{\R^d}
		\frac{\left|\sin\left(\sqrt{\Psi(\xi)}\,
		(r+\epsilon) \right) - \sin\left(
		\sqrt{\Psi(\xi)}\, r\right) \right|^2}{
		\Psi(\xi)} |\hat\phi(\xi)|^2\, d\xi
		\le 8\,\sE ( \epsilon^2;\phi ),
	\end{equation}
	whence $T_1$ is at most $t$ times the right-hand side
	of the preceding. This and \eqref{estimate:T2} together
	yield the proof.
\end{proof}

\begin{proof}[Proof of Theorem \ref{th:cont:t}]
We start with the weak solution $H$ to the stochastic heat
	equation. Throughout, $\phi\in L^2_\L (\R^d)$ is held fixed.
	
	If the integrability condition \eqref{eq:cont:t} holds,
	then according to Lemma \ref{pr:mod-in-t:linear}, for all $T>0$,
	\begin{equation}
		\sup_{\substack{t,\epsilon\ge 0:\\
		0\le t\le t+\epsilon\le T}} \frac{
		\E\left(\left|  H (t+\epsilon\,,\phi)-
		 H (t\,,\phi)\right|^2\right)}{\sE(\epsilon\,;\phi)}
		\le 4 e^{2T}+1<\infty.
	\end{equation}
	Since $\epsilon\mapsto\sE(\epsilon\,;\phi)$ is nondecreasing,
	a direct application of Gaussian-process theory implies
	that $\{ H (t\,,\phi)\}_{t\in[0,T]}$ has a continuous
	modification provided that
	\begin{equation}\label{eq:E:pre}
		\int_{0^+} \frac{\sqrt{\sE(\epsilon\,;\phi)+\F(\epsilon\,;\phi)}}{%
		\epsilon\sqrt{\log(1/\epsilon)}}
		\,d\epsilon<\infty.
	\end{equation}
	See Lemma 6.4.6 of Marcus and Rosen
	\ycite{MarcusRosen}*{p.\ 275}.
	A standard measure-theoretic argument then applies to prove
	that $t\mapsto  H (t\,,\phi)$ has a continuous modification.
	A similar argument works for the weak solution $ W $ to the stochastic
	wave equation \eqref{wave:L}, but we appeal to 
	Proposition \ref{pr:mod-in-t:linear}
	in place of \ref{pr:wave:t-est}.
	Thus, the first portion of our proof
	will be completed, once we prove that condition \eqref{eq:cont:t}
	implies \eqref{eq:E:pre}.
	
	Let us write 
	\begin{equation}\label{cond}\begin{split}
		\int_{0^+} \frac{\sqrt{\sE(\epsilon\,;\phi)+\F(\epsilon\,;\phi)}}{%
			\epsilon\sqrt{\log(1/\epsilon)}}
			\,d\epsilon&\le \int_{0^+} \frac{\sqrt{\sE(\epsilon\,;\phi)}}{%
			\epsilon\sqrt{\log(1/\epsilon)}}\,d\epsilon+
			\int_{0^+} \frac{\sqrt{\F(\epsilon\,;\phi)}}{%
			\epsilon\sqrt{\log(1/\epsilon)}}\,d\epsilon\\
		&:= I_1+I_2.
	\end{split}\end{equation}
	Let us consider $I_1$ first.  We multiply and divide
	the integrand of $I_1$ by the square root of 
	$g(1/\epsilon)$, and then and apply the Cauchy--Schwarz inequality
	to obtain the following:
	\begin{equation}\begin{split}
		I_1&\le \left(\int_{0^+}\frac{\sE(\epsilon\,; \phi)
			g(1/\epsilon)}{\epsilon}\,d\epsilon\right)^{1/2}
			\left(\int_{0^+} \frac{1}{\epsilon \log(1/\epsilon) 
			g(1/\epsilon)}\,d\epsilon\right)^{1/2} \\
	&= \text{const}\cdot \left(\int_{0^+}
			\frac{\sE(\epsilon\,; \phi) g(1/\epsilon)}{\epsilon}
			\,d\epsilon\right)^{1/2}.
	\end{split}\end{equation}
	Note that 
	\begin{equation}\begin{split}
		&\int_0^{1/e} \frac{\sE(\epsilon\,; \phi)
			g(1/\epsilon)}{\epsilon}\,d\epsilon= 
			\frac{1}{(2\pi)^d}\int_{R^d}|\hat{\phi}(\xi)|^2\,d\xi
			\left( \int_0^{1/e}\frac{g(1/\epsilon)}{%
			1+\epsilon \Re \Psi(\xi)}\,d\epsilon \right)\\
		&\hskip1.4in=\frac{1}{(2\pi)^d}\int_{\Re \Psi \geq e}
			|\hat{\phi}(\xi)|^2\,d\xi(\,\cdots)+\frac{1}{(2\pi)^d}
			\int_{\Re \Psi < e}|\hat{\phi}(\xi)|^2\,d\xi(\,\cdots)\\
		&\hskip1.4in:= \frac{I_3+ I_4}{(2\pi)^d},
	\end{split}\end{equation}
	with the notation being clear enough. We will look at $I_3$ and $I_4$ separately. 
	We begin with $I_4$ first.  
	\begin{equation}\label{Ine:1}\begin{split}
		I_4&\le \int_{\Re \Psi<e}\int_0^{1/e}
			\frac{|\hat{\phi}(\xi)|^2g(1/\epsilon)}{%
			1+\epsilon \Re \Psi(\xi)}\,d\epsilon\,d\xi\\
		&\le\int_{\Re \Psi<e}|\hat{\phi}(\xi)|^2\,d\xi
			\cdot \int_0^{1/e}{g(1/\epsilon)}\,d\epsilon\\
		&\le \text{\rm const}\cdot \sE(1\,;\phi),
	\end{split}\end{equation}
	where we have used Proposition \ref{prop:slow} 
	to obtain the last inequality.  In order to bound $I_3$
	we let $N>e$ and look at the following:
	\begin{equation}\label{Ine:2}\begin{split}
		I_3&=\int_0^{1/e}\frac{g(1/\epsilon)}{1+\epsilon N}\,d\epsilon\\
		&=\int_0^{1/N}\frac{g(1/\epsilon)}{1+\epsilon N}
			\,d\epsilon+\int_{1/N}^{1/e}\frac{g(1/\epsilon)}{1+\epsilon N}
			\,d\epsilon\\
		&=I_5+I_6.
	\end{split}\end{equation}
	On one hand, some calculus shows that 
	\begin{equation}\label{Ine:3}
		I_6\le\frac{\log Ng(N)}{N}.
	\end{equation}
	On the other hand, 
	we can change variables and appeal to Proposition \ref{prop:slow} 
	and deduce that
	\begin{equation}\begin{split}
		I_5&\leq \int_N^\infty \frac{g(u)}{u^2}\,du\\
		&\sim \frac{g(N)}{N}\qquad\text{as $N\to\infty$.}
	\end{split}\end{equation}
	We obtain the following upon setting $N:=\Re \Psi(\xi)$ 
	and combining inequalities \eqref{Ine:1}--\eqref{Ine:3} above:
	\begin{equation}\label{eq: cond E}
		\text{\rm \eqref{eq:cont:t}}\quad
		\Longrightarrow\quad
		I_1<\infty.
	\end{equation}
	We now look at $I_2$.  Let us recall the definition of
	$\F(\epsilon\,;\phi)$ and write 
	\begin{equation}\label{Ine:4}\begin{split}
		\F(\epsilon\,;\phi)&=\frac{1}{(2\pi)^d}\left[
			\int_{|\Psi|\leq 1/e}(\,\cdots)\,d\xi+
			\int_{1/e\leq|\Psi|\leq 1/\epsilon}
			(\,\cdots)\,d\xi +\int_{|\Psi|> 1/\epsilon}
			(\,\cdots)\,d\xi \right]\\
		&=\text{\rm const}\cdot[I_7+I_8+I_9].
	\end{split}\end{equation}
	Because $I_7\leq \epsilon^2\sE(1\,;\phi)$ whenever $\epsilon < 1/e$, 
	\begin{equation}\label{Ine:5}
		\int_0^{1/e}\frac{\sqrt{I_7}}{\epsilon \sqrt{\log(1/\epsilon)}}
		\,d\epsilon \leq \text{\rm const}\cdot \sqrt{\sE(1\,;\phi)}.
	\end{equation}
	An application of Cauchy--Schwarz inequality yields
	\begin{equation}\label{Ine:6}\begin{split}
		\int_0^{1/e}\frac{\sqrt{I_8}}{\epsilon
			\sqrt{\log(1/\epsilon)}}\,d\epsilon&=
			\int_0^{1/e}\left(\int_{1/e<|\Psi|\leq 1/\epsilon}
			\frac{|\Psi(\xi)|^2|\hat{\phi}(\xi)|^2}{1+\Re \Psi(\xi)}
			\,d\xi\right)^{1/2} \frac{\epsilon}{\sqrt{\log(1/\epsilon)}}\,d\epsilon\\
		&\le \text{\rm const}\cdot \left(\int_{\R^d}
			\int_0^{1/\Psi(\xi)}\frac{|\Psi(\xi)|^2|\hat{\phi}(\xi)|^2}{%
			1+\Re \Psi(\xi)}\frac{\epsilon^2}{%
			\log(1/\epsilon)}\,d\epsilon \,d\xi\right)^{1/2}\\
		&\le \text{\rm const}\cdot \sqrt{\sE(1\,;\phi)}.
	\end{split}\end{equation}
	We multiply and divide the integrand below by 
	the square root of $g(1/\epsilon)$ and apply
	the Cauchy--Schwarz inequality in order to obtain
	\begin{equation}\label{Ine:7}\begin{split}
		\int_0^{1/e}\frac{\sqrt{I_9}}{\epsilon \sqrt{\log(1/\epsilon)}}
			\,d\epsilon&\le \left(\int_0^{1/e}\int_{|\Psi|>1/\epsilon}
			\frac{|\hat{\phi}(\xi)|^2g(1/\epsilon)}{%
			\epsilon(1+\Re \Psi(\xi))}\,d\xi \,d\epsilon \right)^{1/2}\\
		&\le\text{\rm const}\cdot \left(\int_{\R^d}
			\frac{g(\Psi	(\xi))}{%
			(1+\Re \Psi(\xi))}|\hat{\phi}(\xi)|^2\,d\xi\right)^{1/2}.
	\end{split}\end{equation}
	In the last inequality we have changed the order of
	integration and used Proposition \ref{prop:slow}.
	Taking into account inequalites (\ref{Ine:3})--(\ref{Ine:7}), we obtain
	\begin{equation}\label{eq: cond F}
		\quad\left(\int_{\R^d}
		\frac{g(1+\Psi(\xi))}{(1+\Re \Psi(\xi))}
		|\hat{\phi}(\xi)|^2\,d\xi\right)^{1/2}<\infty\quad
		\Longrightarrow\quad
		I_2<\infty.
	\end{equation}
	Inequalities $\eqref{eq: cond E}$ and $\eqref{eq: cond F}$
	together with $\eqref{cond}$ imply that 
	\begin{equation}  
		\int_{0^+} \frac{\sqrt{\sE(\epsilon \,;\phi)+
		\F(\epsilon\,;\phi)}}{\epsilon \sqrt{\log(1/\epsilon)}}
		\, d\epsilon<\infty,
	\end{equation}
	provided that \eqref{eq:cont:t} holds.
	This proves the first assertion of the theorem.

	Suppose $\gamma>0$, where
	\begin{equation}
		\gamma:=\underline{\text{\rm ind}}\,\sE(\bullet\,;\phi),
	\end{equation}
	for brevity. Then by definition,
	$\sE(\epsilon\,;\phi)\le \epsilon^{\gamma+o(1)}$ as $\epsilon\downarrow 0$.
	Another standard result from Gaussian analysis, used in conjunction
	with Proposition \ref{pr:mod-in-t:linear}
	proves that $ H $ has
	a H\"older-continuous modification with H\"older exponent
	$\le \gamma/2$ \cite{MarcusRosen}*{Theorem 7.2.1, p.\ 298}.
	The proof for $ W $ is analogous.
	
	For the remainder of the proof we consider only the weak solution
	$ H $ to the stochastic heat equation \eqref{heat:L}, and
	write $H_t :=  H (t\,,\phi)$
	for typographical ease. If $\gamma>0$,
	then, Proposition \ref{pr:mod-in-t:linear} and elementary
	properties of normal laws together imply that
	\begin{equation}
		\inf_{t\ge 0} \left\| H_{t+\epsilon}-H_t \right\|_{L^2(\P)}
		\ge \epsilon^{\gamma+o(1)}\quad\text{for infinitely-many
		$\epsilon\downarrow 0$}.
	\end{equation}
	Consequently, for all $\delta\in(0\,,1)$, $q>\gamma/2$, $T>S>0$,
	and $t\in[S\,,T]$---all fixed---the following holds for infinitely-many
	values of $\epsilon\downarrow 0$:
	\begin{equation}\begin{split}
		\P\left\{ \sup_{r\in[S,T]}\left| H_{r+\epsilon} - H_r \right|
			\ge \epsilon^q \right\} &\ge 
			\P\left\{ \left| H_{t+\epsilon} - H_t \right|
			\ge \delta \left\| H_{t+\epsilon} -
			H_t \right\|_{L^2(\P)}\right\}\\
		&= 1- \sqrt{\frac{2}{\pi}}\int_0^\delta \exp(-x^2/2)\, dx.
	\end{split}\end{equation}
	The inequality $\exp(-x^2/2)\le1$ then implies
	that
	\begin{equation}
		\P\left\{
		\sup_{r\in[S,T]} | H_{r+\epsilon}-H_r |\ge (1+o(1)) \epsilon^q\
		\text{for infinitely many $\epsilon\downarrow 0$}
		\right\} \ge 1-\sqrt{\frac{2}{\pi}}\,\delta.
	\end{equation}
	Since $q>\gamma/2$ is arbitrary, we can enlarge it arbitrarily
	to infer that
	\begin{equation}
		\P\left\{ \limsup_{\epsilon\downarrow 0} \sup_{r\in[S,T]}\frac{
		\left| H_{r+\epsilon}-H_r\right|}{\epsilon^q}=\infty
		\right\}\ge 1-\sqrt{\frac{2}{\pi}}\,\delta.
	\end{equation}
	Let $\delta\downarrow 0$ to find that 
	\begin{equation}\label{eq:H:LB}
		\limsup_{\epsilon\downarrow 0} \sup_{r\in[S,T]}\frac{
		\left| H_{r+\epsilon}-H_r\right|}{\epsilon^q}=\infty
		\qquad\text{a.s.}
	\end{equation}
	This proves that any $q>\gamma/2$ is an almost-sure lower bound
	for the critical H\"older exponent of $H$. Moreover,
	in the case that $\gamma=0$, we find that
	\eqref{eq:H:LB} holds a.s.\ for all $q>0$. Thus, it
	follows that with probability one,
	$H$ has no H\"older-continuous modification in that case.
	The proof is now complete.
\end{proof}

\section{\bf Spatial and joint continuity: Proofs of Theorems \ref{th:Holder}
	and \ref{th:joint:holder}}\label{sec6}

\begin{proof}[Proof of Theorem \ref{th:Holder}]
		We begin by proving the portion of
		Theorem \ref{th:Holder} that relates to the stochastic
		heat equation and its random-field solution
		$\{ H (t\,,x);\, t\ge 0,\, x\in\R\}$; 
		Theorem \ref{th:exist} guarantees the existence
		of the latter process. 

		Throughout we can---and will---assume without loss of generality
		that $x\mapsto  H (t\,,x)$ is continuous, and hence
		so is the local time of $\bar{X}$ in its spatial variable.
		As we saw, during the course of the proof of Theorem \ref{th:cont},
		this automatically implies the condition \eqref{eq:LT:E}, which
		we are free to assume henceforth.

		According to Proposition \ref{pr:heat:embed},
		specifically Corollary \ref{cor:heat:embed},
		\begin{equation}
			\tfrac{1}{3}\sE(t\,;\delta_x-\delta_y)\le
			\E\left( \left|  H (t\,,x)- H (t\,,y) \right|^2\right)
			\le 4 \sE(t\,;\delta_x-\delta_y).
		\end{equation}
		Since $|\hat\delta_x(\xi)-\hat\delta_y(\xi)|^2=2[1-\cos(\xi(x-y))]$,
		it follows from this and \eqref{eq:Omega} that
		\begin{equation}\label{eq:heat:hHh}
			\tfrac{2}{3} h (|x-y|)\le
			\E\left( \left|  H (t\,,x)- H (t\,,y) \right|^2\right)
			\le 8 h (|x-y|).
		\end{equation}
		This implies that the critical
		H\"older exponent of $z\mapsto H (t\,,z)$ is
		almost surely equal to one-half of the following quantity:
		\begin{equation}\label{eq:h}
			\underline{\text{\rm ind}}\,  h  := \liminf_{\epsilon\downarrow 0}
			\frac{\log h (\epsilon)}{\log\epsilon}.
		\end{equation}
		We will not prove this here, since it is very similar to the proof
		of temporal H\"older continuity (Theorem \ref{th:cont:t}). Consequently,
		\begin{equation}\label{eq:Holder:iff:Gauss}
			z\mapsto  H (t\,,z)\text{ has a H\"older-continuous
			modification iff }\underline{\text{ind}}\,  h >0.
		\end{equation}
		Among other things, this implies the equivalence of parts
		(1) and (2) of the theorem.

		Let $ Z (t\,,x)$ denote the local time of $\bar X$ at spatial value
		$x$ at time $t\ge 0$. We prove that (1), (2), and (3) are equivalent
		by proving that \eqref{eq:Holder:iff:Gauss} continues to hold
		when $ H $ is replaced by $ Z $. Fortunately, this can be
		read off the work of Barlow \ycite{Barlow}.
		We explain the details briefly.
		Because \eqref{eq:LT:E} holds,
		Theorem 5.3 of Barlow \ycite{Barlow} implies
		that there exists a finite constant $c>0$
		such that for all $t\ge 0$ and finite intervals $I\subset\R$,
		\begin{equation}
			\lim_{\delta\downarrow 0}\sup_{\substack{a,b\in I\\
			|a-b|<\delta}} \frac{| Z (t\,,a)- Z (t\,,b)|}{
			\sqrt{ h (|a-b|) \log(1/|b-a|)}}\ge c
			\left(\sup_{x\in I}
			 Z (s\,,x)\right)^{1/2}\qquad\text{a.s.}
		\end{equation}
		Moreover, we can choose the null set to be independent of all intervals
		$I\subset\R$ with rational endpoints. 
		In fact, we can replace $I$ by $\R$, since $a\mapsto Z (t\,,a)$
		is supported by the closure of the
		range of the process $\bar{X}$ up to time $t$,
		and the latter range is a.s.\ bounded since $\bar{X}$ is cadlag.

		By their very definition local times satisfy
		$\int_{-\infty}^\infty  Z (s\,,x)\, dx=s$ a.s.
		Thus, $\sup_{x\in\R}
		 Z (t\,,x)>0$ a.s., whence it follows
		that for all $q>\frac12\,\underline{\text{\rm ind}}\,
		 h $, 
		\begin{equation}
			\lim_{\delta\downarrow 0}\sup_{ |a-b|<\delta}
			\frac{| Z (t\,,a)- Z (t\,,b)|}{
			|a-b|^q}=\infty\qquad\text{a.s.}
		\end{equation}
		That is, there is no H\"older-continuous modification of
		$a\mapsto Z (t\,,a)$ of order $>\frac12\,\underline{\text{ind}}\,
		 h $. In particular, if $\underline{\text{ind}}\,
		 h =0$, then $a\mapsto Z (t\,,a)$ 
		does not have a H\"older-continuous 
		modification.

		Define $d(a\,,b):=\sqrt{ h (|a-b|)}$; it is easy to see
		that $d$ is a pseudo-metric on $\R$.
		According to Bass and Khoshnevisan \ycite{BK},
		\begin{equation}
			\limsup_{\delta\downarrow 0}\sup_{d(a,b)<\delta}
			\frac{| Z (t\,,a)- Z (t\,,b)|}{
			\int_0^{d(a,b)} (\log N(u))^{1/2}\,du} 
			\le 2\left(\sup_x Z (t\,,x)\right)^{1/2}\qquad\text{a.s.,}
		\end{equation}
		where $N(u)$ denotes the smallest number of 
		$d$-balls of radius $\le u$ needed to cover $[-1\,,1]$.
		This sharpened an earlier result of 
		Barlow \ycite{Barlow:85}*{Theorem 1.1}.
		Furthermore,
		$\sup_x Z (t\,,x)<\infty$ a.s.\
		\cite{BK}*{Theorem 3.1}.
		Consequently,
		\begin{equation}
			\sup_{d(a,b)<\delta}
			| Z (t\,,a)- Z (t\,,b)|=
			O\left(\int_0^\delta (\log N(u))^{1/2}\, du\right)\qquad
			\text{as $\delta\downarrow 0$}.
		\end{equation}
		Since $ h $ is increasing,
		\begin{equation}\label{eq:HC:L:UB}
			\sup_{|a-b|<\delta}
			| Z (t\,,a)- Z (t\,,b)|=
			O\left(\int_0^{ h ^{-1}(\delta)} (\log N(u))^{1/2}\, du\right)\qquad
			\text{as $\delta\downarrow 0$}.
		\end{equation}

		According to equations (6.128) and (6.130) of Marcus and Rosen
		\ycite{MarcusRosen}*{Lemma 6.4.1, p.\ 271}, there exists a finite constant
		$c>0$ such that for all $u>0$ small,
		\begin{equation}\begin{split}
			N(u) &\le \frac{c}{\ell_2 \left\{ (x\,,y)\in[-1\,,1]^2:\
				h (|x-y|)<u/4\right\}}\\
			&\le \frac{\text{const}}{ h ^{-1}(u/4)},
		\end{split}\end{equation}
		where $\ell_2(A)$ denotes the Lebesgue measure of $A\subset\R^2$.
		[Specifically, we apply Lemma 6.4.1 of that reference with
		their $K:=[-1\,,1]$ and their $\mu_4:=c$.] This and \eqref{eq:HC:L:UB}
		together imply that with probability one the following is valid:
		As $\delta\downarrow 0$,
		\begin{equation}\begin{split}
			\sup_{|a-b|<\delta}
				| Z (t\,,a)- Z (t\,,b)| &=
				O\left(\int_0^\delta |\log u|^{1/2}\,  h (du)\right)\\
			&=O\left( h (\delta)|\log(1/\delta)|^{1/2}\right)+
				O\left(\int_0^\delta \frac{ h (u)}{u|\log u|^{1/2}}\, du\right).
		\end{split}\end{equation}
		The last line follows from integration by parts.
		If $\gamma:=\underline{\text{\rm ind}}\, h >0$, then 
		$ h (\delta)=o(\delta^q)$ for
		all fixed choices of $q\in(0\,,\gamma/2)$. It follows that
		$a\mapsto  Z (t\,,a)$ is H\"older continuous of
		any order $<\gamma/2$. Among other things, this implies
		\eqref{eq:Holder:iff:Gauss}  with $ Z $
		replacing $ H $, whence it follows that (1)--(3)
		of the theorem are equivalent.

		The hyperbolic portion of the theorem is proved similarly,
		but we use Proposition \ref{pr:wave:embed} in place
		of Proposition \ref{pr:heat:embed} everywhere.
\end{proof}

\begin{proof}[Proof of Theorem \ref{th:joint:holder}]
	Since $\beta''$ is positive, $\Re\Psi(\xi)\to\infty$
	as $|\xi|\to\infty$. Therefore,
	for all $\epsilon,\vartheta>0$ and $x\in\R$,
	\begin{equation}\begin{split}
		\sE(\epsilon\,;\delta_x) &=\frac{1}{(2\pi)^d}
			\int_{\R^d}\frac{d\xi}{(1/\epsilon)+\Re\Psi(\xi)}\\
		&\leq \text{const}\cdot\epsilon +
			\int_{\{\Re\Psi>\vartheta\}}
			\frac{d\xi}{(1/\epsilon)+\Re\Psi(\xi)}.
	\end{split}\end{equation}
          If $\xi$ is sufficiently large, then
	for all $\gamma\in(d\,,\beta'')$,
	we can find a constant $C_\gamma\in(0\,,\infty)$ such that
	$\Re\Psi(\xi)\ge C_\gamma|\xi|^\gamma$ for all
	$\xi\in\{\Re\Psi>\vartheta\}$. Consequently,
	\begin{equation}\begin{split}
		\sE(\epsilon\,;\delta_x) &= O\left(
			\epsilon +
			\int_{\{\Re\Psi>\vartheta\}}\frac{d\xi}{(1/\epsilon)+
			|\xi|^\gamma }
			\right)\\
		&=O\left( \epsilon^{1-(d/\gamma)}\right),
	\end{split}\end{equation}
	as $\epsilon\downarrow 0$. Thus,
	$\underline{\text{\rm ind}}\, \sE(\bullet\,;\delta_x)
	\le 1-(d/\beta'')$, where this index was introduced
	in \eqref{eq:Hcont:t}. A similar calculation shows that
	$\underline{\text{\rm ind}}\, h \le \beta''-d$;
	confer with \eqref{eq:h} for the definition of this quantity. Thus,
	for all fixed $T>0$,
	Proposition \ref{pr:mod-in-t:linear} and \eqref{eq:heat:hHh}
	together prove the following: For all $\tau<(\beta''-d)/\beta''$,
	$\zeta<\beta''-d$, $x,y\in\R$, and $s,t\in[0\,,T]$,
	\begin{equation}\label{joint:holder}
		\E\left(\left|  H (t\,,x)- H (s\,,y)
		\right|^2\right)\le \text{const}\cdot\left(
		|s-t|^\tau + |x-y|^\zeta\right).
	\end{equation}
	A two-dimensional version of Kolmogorov's continuity theorem
	finishes the proof; see the proof of
	Theorem (2.1) of Revuz and Yor \ycite{RY}*{p.\ 25}.
	
	The proof, in the case of the stochastic
	wave equation, is similar, but we use Propositions
	\ref{pr:wave:embed} and 
	\ref{pr:mod-in-t:linear} instead of Propositions 
	\ref{pr:heat:embed} and \ref{pr:wave:t-est}, respectively.
\end{proof}

\section{\bf Heat equation via generators of Markov processes}\label{sec7}
We now consider briefly the stochastic heat equation,
where the spatial movement is governed by the generator $\L$
of a [weakly] Markov process $X:=\{X_t\}_{t\ge 0}$ that takes values in
a locally compact separable metric space $F$.
We assume further that $X$ admits a symmetrizing measure
$m$ that is Radon and positive. 
Let us emphasize that
$m$ satisfies $(P_tf,g)=(f,P_t g)$ for all
$t\ge 0$ and $f,g\in L^2(m)$, where $\{P_t\}_{t\ge 0}$ denotes
the transition operators of $X$, and $(\phi_1\,,\phi_2):=\int 
\phi_1\phi_2\,dm$
for all $\phi_1,\phi_2\in L^2(m)$.

\subsection{The general problem}
Consider the stochastic heat equation
\begin{equation}\label{heat1}\left|\begin{split}
	&\partial_t u (t\,,x) = (\L  u)(t\,,x)
		+ \dot{w}(t\,,x),\\
	&u(0\,,x) = 0,
\end{split}\right.\end{equation}
valid for all $t\ge 0$ and $x\in F$.
Here, the underlying noise $w$ in \eqref{heat1}
is a Gaussian martingale measure
on $\R_+\times F$ in the sense of Walsh \ycite{Walsh}:
$w$ is defined on the filtered probability space
$( \Omega \,, \F,\{ \F_t\}_{t\ge 0},\P)$ and
$w_t(\phi) := \int_0^t\int_F \phi(s\,,x)\, w(dx\,ds)$
defines an $\{\F_t\}_{t\ge 0}$-martingale for
$\phi\in L^2(ds\times m)$;
$w$ can be characterized by the covariance functional
for the corresponding Wiener integrals:
\begin{equation}\label{def2}
	\E\left(\int f\, dw\cdot \int g\, dw \right)=\int_0^\infty
	\int_F \int_F f(s\,,x) g(s\,,y)
	\,m(dx)\,m(dy)\,ds,
\end{equation}
for all $f,g\in L^2(ds\times m)$.

We can follow the description of Walsh \ycite{Walsh}
and write the weak form of
equation $\eqref{heat1}$ as follows:
For all $\phi\in L^2(m)$ and $t\ge 0$, 
\begin{equation}\label{wsol}
	u(t\,,\phi)=\int_0^t\int_F
	(P_{t-s} \phi) (x)\,w(dx\,ds).
\end{equation}

\begin{lemma}
	The integral defined by \eqref{wsol} is well defined
	for all $\phi\in L^2(m)$.
\end{lemma}

\begin{proof}
	We follow closely the proof of Proposition \ref{pr:heat:exist},
	and apply the fact that the semigroup $\{P_s\}_{s\ge 0}$ 
	is a contraction on $L^2(m)$.
\end{proof}

Let $Z$ denote the
occupation measure of $X$; consult \eqref{OM}.
The following is the key result of this section. It 
identifies an abstract Hilbertian quasi-isometry between
the occupation-measure $L^2$-norm of $X$ and a similar
norm for the solution to
the stochastic heat equation \eqref{heat1} for $\L$.

\begin{theorem}\label{th:LT:u}
	If $u$ denotes the weak solution to
	\eqref{heat1}, then for all $\phi\in L^2(m)$ and $t\ge 0$,
	\begin{equation}
		\tfrac{1}{8}\,t\E\left(\left| u(t\,,\phi)\right|^2\right)
		\le \E_m\left(\left|  Z(t\,,\phi) \right|^2\right)
		\le 4t  \E\left(\left|u(t\,,\phi) \right|^2\right).
	\end{equation}
\end{theorem}

As usual, $\E_m$ refers to the expectation operator for the 
process $X$, started according to the measure $m$.

The preceding theorem follows from the next formula.

\begin{proposition}\label{pr:L2:1}
	If $u$ denotes the weak solution to
	\eqref{heat1}, then for all  $\phi\in L^2(m)$ and $t\ge 0$,
	\begin{equation}\label{eq:L2:1}
		\E_m\left(\left| Z(t\,,\phi) \right|^2\right)
		=4 \int_0^t \E\left(\left|u(s/2\,,\phi) \right|^2\right)\,ds.
	\end{equation}
\end{proposition}

\begin{proof}[Proof of Proposition \ref{pr:L2:1}]
	Since $m$ is a symmetrizing measure for $X$, the $\P_m$-law
	of $X_u$ is $m$ for all $u\ge 0$. By the Markov
	property and Tonelli's theorem,
	\begin{equation}
		\E_m\left(\left|  Z(t\,,\phi)\right|^2\right)=
		2\int_0^t\int_u^t \left( P_{v-u}\phi\,,\phi\right)\, dv\,du.
	\end{equation}
	We computing the Laplace transform of both sides, viz.,
	\begin{equation}\begin{split}
		\int_0^\infty e^{-\lambda t}
			\E_m\left(\left|  Z(t\,,\phi)\right|^2\right)
			\, dt &= 2\int_0^\infty
			\int_0^t \int_u^t  e^{-\lambda t}
			\left( P_{v-u}\phi\,,\phi \right)\, dv\,du\,dt\\
		&=2\int_0^\infty\int_u^\infty\left(\int_v^\infty e^{-\lambda t}
			\, dt\right) \left( P_{v-u}\phi\,,\phi \right)\, dv\, du\\
		&=\frac{2}{\lambda^2}\int_0^\infty e^{-\lambda s}
			\left( P_s\phi\,,\phi\right)\, ds.
	\end{split}\end{equation}
	The exchange of the integrals is justified because
	$(P_r\phi\,,\phi)=\|P_{r/2}\phi\|_{L^2(m)}^2$ is positive and finite.
	Because $\phi\in L^2(m)$, Fubini's theorem implies that for all $\lambda>0$,
	\begin{equation}
		\int_0^\infty e^{-\lambda t}
		\E_m\left(\left|  Z(t\,,\phi)\right|^2\right)
		\, dt = \frac{2}{\lambda^2}\left( R_\lambda\phi\,,\phi\right),
	\end{equation}
	where $R_\lambda:=\int_0^\infty \exp(-\lambda s)
	P_s\, ds$ defines
	the resolvent of $\{P_t\}_{t\ge 0}$.
	
	Let $T_\lambda$ denote an independent mean-$(1/\lambda)$ exponential
	holding time. The preceding display can be rewritten as follows:
	\begin{equation}\label{eq:Em:L:R}
		\E_m\left(\left|  Z(T_\lambda\,,\phi)\right|^2\right)
		=\frac{2}{\lambda}\left(R_\lambda\phi\,,\phi\right)
		\qquad\text{for all $\lambda>0$}.
	\end{equation}
	
	Next we consider the weak solution $u$ to the stochastic
	heat equation \eqref{heat1} by first observing that
	$\E ( | u(t\,,\phi) |^2 )=
	\int_0^t\|P_s\phi\|_{L^2(m)}^2\,ds$. It follows from
	this, and successive applications of Tonelli's theorem, that for all $\beta>0$,
	\begin{equation}\begin{split}
		\E\left(\left| u(T_\beta\,,\phi)\right|^2\right)
			&= \int_0^\infty \beta e^{-\beta t} \int_0^t
			\|P_s\phi\|_{L^2(m)}^2\, ds\, dt\\
		&=\int_0^\infty e^{-\beta s}\|P_s\phi\|_{L^2(m)}^2\, ds.
	\end{split}\end{equation}
	Because $\|P_s\phi\|_{L^2(m)}^2 =
	(P_{2s}\phi\,,\phi)$, we may apply Fubini's theorem once more,
	and select $\beta:=2\lambda$, to find that
	\begin{equation}\label{eq:E:L:u}\begin{split}
		\E\left(\left| u(T_{2\lambda} \,,\phi)\right|^2\right)
			&= \int_0^\infty e^{-2\lambda s} \left( P_{2s}\phi\,,\phi
			\right)\,ds\\
		&= \left( \phi\,,\int_0^\infty e^{-2\lambda s} P_{2s}\phi\, ds\right)\\
		&=\frac12\left(R_\lambda\phi\,,\phi\right).
	\end{split}\end{equation}
	The condition of square integrability for $\phi$ justifies
	the appeal to Fubini's theorem.
	We can compare \eqref{eq:Em:L:R} and \eqref{eq:E:L:u} to find that
	\begin{equation}
		\E_m\left(\left|  Z(T_\lambda\,,\phi)\right|^2\right)
		=\frac{4}{\lambda}\E\left(\left| u(T_{2\lambda}\,,\phi)
		\right|^2\right)\qquad\text{for all $\lambda>0$}.
	\end{equation}
	Define
	$q(t):=\E_m(| Z(t\,,\phi)|^2)$,
	$\rho(t):=\E(|u(t/2\,,\phi)|^2)$
	and $\1(t):=1$ for all $t\ge 0$. The preceding shows
	that the Laplace transform of $q$ is equal to $4$ times
	the product of the respective Laplace transforms of $\rho$
	and $\1$. Thus, we can invert to find that
	$q=4\rho*\1$, which is another way to state the theorem.
\end{proof}

\begin{proof}[Proof of Theorem \ref{th:LT:u}]
	Let us choose and fix a measurable function $\phi:F\to\R$ such that 
	$|\phi|\in L^2(m)$. The defining isometry for Wiener integrals
	yields the following identity, where both sides are convergent:
	$\E ( | u(t\,,\phi) |^2 ) 
	= \int_0^t\|P_s \phi\|^2_{L^2(m)}\,ds$.
	Proposition \ref{pr:L2:1} implies that
	\begin{equation}\label{eq:7.17}
		2t\E\left( \left| u(t/4\,,\phi)\right|^2\right)\le
		\E_m\left(\left|  Z(t\,,\phi)\right|^2\right)
		\le 4t\E\left( \left| u(t/2\,,\phi) \right|^2\right).
	\end{equation}
	[For the lower bound, we use the bound
	$\int_0^t\E(|u(s/2\,,\phi)|^2)\, ds\ge
	\int_{t/2}^t\E(|u(s/2\,,\phi)|^2)\, ds$.]	
	By monotonicity, $\E(|u(t/2\,,\phi)|^2)\le \E(|u(t\,,\phi)|^2)$,
	whence follows the announced
	upper bound for $\E_m(| Z(t\,,\phi)|^2)$. 
	
	In order to prove the other bound we first write
	\begin{equation}
		 Z(t\,,\phi) =
		 Z(t/2\,,\phi)
		+ Z(t/2\,,\phi) \circ\theta_{t/2},
	\end{equation}
	where $\{\theta_s\}_{s\ge 0}$ denotes the collection of all shifts
	on the paths of $X$. We care only about distributional properties.
	Therefore, by working on an appropriate probability space,
	we can always insure that these shifts can be constructed;
	see Blumenthal and Getoor \ycite{BG:68}.
	 
	We apply the Markov property at time $t/2$.
	Since $\P_m\circ X_{t/2}^{-1}=m$, it follows that
	\begin{equation}
		\E_m \left( \left| Z (t/2\,,\phi )
		\circ\theta_{t/2} \right|^2
		\right)=
		\E_m \left( \left| Z (t/2\,,\phi )
		\right|^2 \circ\theta_{t/2} \right)=
		\E_m\left(\left| Z (t/2\,,\phi )
		\right|^2\right),
	\end{equation}
	and hence
	$\E_m(| Z(t\,,\phi)|^2)\le 4\E_m(| Z(t/2\,,\phi)|^2)$.
	Consequently,
	\begin{equation}
		\E_m\left(\left|
		 Z(t\,,\phi) \right|^2\right)
		\le 16\E_m \left(\left| Z(t/4\,,\phi) \right|^2\right)
		\qquad\text{for all $t\ge 0$}.
	\end{equation}
	This and the first inequality of \eqref{eq:7.17} together imply the remaining
	bound in the statement of the theorem.
\end{proof}

\subsection{The stochastic heat equation in dimension $2-\epsilon$}

We now specialize the setup of the preceding subsection to
produce an interesting family of examples:
We suppose that $F$ is a
locally compact subset of $\R^d$ for some integer $d\ge 1$,
and $m$ is a positive Radon measure on $F$, as before.
Let $\{R_\lambda\}_{\lambda>0}$ denote the resolvent of $X$,
and suppose that $X$ has jointly continuous and uniformly
bounded resolvent densities
$\{r_\lambda\}_{\lambda>0}$.
In particular, $r_\lambda(x\,,y)\ge 0$ for all $x,y\in F$ and
\begin{equation}
	(R_\lambda f)(x) =\int_F r_\lambda(x\,,y)f(y)\, m(dy),
\end{equation} 
for all measurable functions $f:F\to\R_+$. Recall that $X$
has local times $\{ Z(t\,,x)\}_{t\ge 0,x\in F}$
if and only if for all measurable functions
$f:F\to\R_+$, and every $t\ge 0$,
\begin{equation}
	 Z(t\,,f)=\int_F  Z(t\,,z)f(z)\, m(dz),
\end{equation}
valid $\P_x$-a.s.\ for all $x\in F$. Choose and fix some
point $a\in F$, and define
\begin{equation}\label{eq:fa}
	f_\epsilon^a (z) := \frac{\1_{B(a,\epsilon)}(z)}{m(B(a\,,\epsilon))}
	\qquad\text{for all $z\in F$ and $\epsilon>0$}.
\end{equation}
Of course, $B(a\,,\epsilon)$ denotes the ball of radius $\epsilon$
about $a$, measured in the natural metric of $F$. Because 
$r_\lambda$ is jointly continuous,
$\lim_{\epsilon\downarrow 0}(R_\lambda f_\epsilon^a)(x)=r_\lambda(x\,,a)$,
uniformly for $x$-compacta. Define $\phi_{\epsilon,\delta}:=f_\epsilon^a-f_\delta^a$,
and observe that $\phi\in L^2(m)$. Furthermore,
\begin{equation}\label{RRR}
	\lim_{\epsilon,\delta\downarrow 0}\left( R_\lambda\phi_{\epsilon,\delta}
	\,,\phi_{\epsilon,\delta}\right) = \lim_{\epsilon,\delta\downarrow 0}
	\left\{ \left( R_\lambda f_\epsilon^a\,,f_\epsilon^a\right)
	-2\left( R_\lambda f_\epsilon^a\,, f_\delta^a \right) 
	+\left( R_\lambda f_\delta^a \,,f_\delta^a \right)\right\}
	=0.
\end{equation}
If $h\in L^2(m)$, then the weak solution $h$ to \eqref{heat1}---where
$\L$ denotes the $L^2$-generator of $X$---satisfies
\begin{equation}\begin{split}
	\E\left( \left| u(t\,,h) \right|^2\right) &=
		\int_0^t \|P_s h\|_{L^2(m)}^2\, ds\\
	&\le e^{2\lambda t}\int_0^\infty e^{-2\lambda s}\|P_sh\|_{L^2(m)}^2\, ds\\
	&=\frac{e^{2\lambda t}}{2} (R_\lambda h\,,h).
\end{split}\end{equation}
Therefore, Theorem \ref{th:LT:u} and \eqref{RRR} together
imply that $\{u(t\,,f_\epsilon^a)\}_{\epsilon>0}$
is a Cauchy sequence in $L^2(\P)$ for all $t\ge 0$.
In other words, we have shown that the stochastic heat equation
\eqref{heat1} has a ``random-field solution.''

\begin{example}\label{example:fractal}
	It is now easy to check, using the heat-kernel estimates of
	Barlow
	\ycite{Barlow:99}*{Theorems 8.1.5 and 8.1.6}, that for all $d \in(0\,,2)$
	there exists a compact ``fractal'' $F\subset\R^2$ of
	Hausdorff dimension $d$ such that Brownian motion
	on $F$ satisfies the bounded/continuous resolvent-density
	properties here. The preceding proves that if we replace
	$\L$ by the Laplacian on $F$, then the stochastic
	heat equation \eqref{heat1} has a ``random-field solution.''
	Specifically, the latter means that  for all $t\ge 0$ and
	$a\in F$,
	$u(t\,,f_\epsilon^a)$ converges in $L^2(\P)$ as $\epsilon\downarrow 0$,
	where $f_\epsilon^a$ is defined in \eqref{eq:fa}. This example
	comes about, because the fractional diffusions of 
	Barlow \ycite{Barlow:99}
	have local times when the dimension $d$ of the fractal on which they live
	satisfies $d<2$. See Barlow
	\ycite{Barlow:99}*{Theorem 3.32}, for instance.
	\qed
\end{example}

\section{\bf A semilinear parabolic problem}\label{sec8}
We consider the semilinear problems that correspond to 
the stochastic heat equation \eqref{heat:L}. 
At this point in the development of SPDEs, we can make
general sense of nonlinear stochastic PDEs only when the linearized SPDE
is sensible. Thus, we assume henceforth that $d=1$. 

We investigate the
semilinear stochastic heat equation. Let $b:\R\to\R$ be a measurable
function, and consider the solution $ H_b$ to the following
SPDE:
\begin{equation}\label{heat:semilinear}\left|\begin{split}
	&\partial_t  H_b(t\,,x) = (\L  H_b) (t\,,x) 
		+ b( H_b(t\,,x)) + \dot{w}(t\,,x),\\
	& H_b(0\,,x)=0,
\end{split}\right.\end{equation}
where $\L$ denotes the generator of the L\'evy process $X$, as before.

Equation \eqref{heat:semilinear} has a chance of making 
sense only if the linearized problem \eqref{heat:L}
has a random-field solution $ H $, in which case
we follow Walsh \ycite{Walsh} and write the solution $ H_b$
as the solution to the following:
\begin{equation}\label{eq:Hb:H:1}
	 H_b(t\,,x) =  H (t\,,x) + \int_0^t\int_{-\infty}^\infty
	b( H_b(s\,,x-y))\, P_{t-s}(dy)\, ds,
\end{equation}
where the measures $\{P_t\}_{t\ge 0}$ are determined
from the semigroup of $X$ by $P_t(E):=(P_t\1_E)(0)$ for all
Borel sets $E\subset\R$. [This is standard notation.]
We will soon see that this random integral equation has a
``good solution'' $ H_b$
under more or less standard conditions on the function $b$.
But first, let us make an observation.

\begin{lemma}\label{lem:heat:tdf}
	If \eqref{heat:L} has a random-field solution, then
	the process $X$ has a jointly measurable transition
	density $\{p_t(x)\}_{t> 0, x\in\R}$ that satisfies
	the following: For all $\eta>0$ there exists a constant
	$C:=C_\eta\in(0\,,\infty)$ such that for all $t>0$,
	\begin{equation}
		\int_0^t \|p_s\|_{L^2(\R)}^2\, ds\le C e^{\eta t}.
	\end{equation}
	Finally, $(t\,,x)\mapsto p_t(x)$ is uniformly continuous
	on $[\epsilon\,,T]\times\R$ for all fixed $\epsilon,T>0$.
\end{lemma}

\begin{proof}
	We can inspect the function $y=
	x\exp(-x)$ to find that $\exp(-x)\le (1+x)^{-1}$ for
	all $x\ge 0$. Consequently,
	\begin{equation}\label{eq:firstsecond}
		\int_{-\infty}^\infty e^{-\epsilon\Re\Psi(\xi)}\, d\xi
		\le \int_{-\infty}^\infty \frac{d\xi}{1+\epsilon
		\Re\Psi(\xi)}.
	\end{equation}
	The second integral, however, has been shown to be equivalent to the
	existence of random-field solutions to \eqref{heat:L};
	see \eqref{cond:hawkes}. It follows that the first integral in
	\eqref{eq:firstsecond} is convergent. We apply the inversion theorem
	to deduce from this that the transition densities of $X$ are
	given by $p_t(x) = (2\pi)^{-1}\int_{-\infty}^\infty
	\exp \{-ix\xi - t\Psi(\xi )\} d\xi$,
	where the integral
	is absolutely convergent for all $t>0$ and $x\in\R$. Among other things,
	this formula implies the uniform continuity of $p_t(x)$ away from
	$t=0$. In addition,
	by Plancherel's theorem, for all $s>0$,
	\begin{equation}\begin{split}
		\|p_s\|_{L^2(\R)}^2 &= \frac{1}{2\pi}\int_{-\infty}^\infty
			\left| e^{-s\Psi(\xi)} \right|^2\,d\xi\\
		&=  \frac{1}{2\pi}\int_{-\infty}^\infty e^{-2s\Re\Psi(\xi)}\, d\xi.
	\end{split}\end{equation}
	Therefore, Lemma \ref{lem:Monotone} 
	and Tonelli's theorem together imply that for all $\lambda>0$,
	\begin{equation}
		\int_0^t \|p_s\|_{L^2(\R)}^2\, ds 
		\le  \frac{e^{2t/\lambda}}{4\pi}\int_{-\infty}^\infty
		\frac{d\xi}{(1/\lambda)+\Re\Psi(\xi)}<\infty.
	\end{equation}
	This completes the proof.
\end{proof}
 
Thanks to the preceding lemma, 
by \eqref{heat:semilinear} we mean a solution to the following:
\begin{equation}\label{eq:Hb:H:2}
	 H_b(t\,,x) =  H (t\,,x) + \int_0^t\int_{-\infty}^\infty
	b( H_b(s\,,x-y))\, p_{t-s}(y)\, dy\, ds.
\end{equation}

For the following we assume that the underlying probability
space $(\Omega\,,\F,\P)$ is complete.

\begin{theorem}\label{th:heat:semilinear}
	Suppose $b$ is bounded and globally Lipschitz,
	and the stochastic heat equation
	\eqref{heat:L} has a random-field solution $H$;
	thus, in particular, $d=1$. Then,
	there exists a modification of $H$, denoted still by $H$,
	and a process $ H_b$ with the
	following properties:
	\begin{enumerate}
		\item $ H_b\in L^p_{\text{\it loc}}(\R_+\times\R)$
			for all $p\in[1\,,\infty)$.
		\item With probability one, \eqref{eq:Hb:H:2} holds for all
			$(t\,,x)\in\R_+\times\R$.
		\item For all $T>0$,
			$J$ is a.s.\ bounded
			and continuous on $[0\,,T]\times\R$, where
			\[
				J(t\,,x) := \int_0^t\int_{-\infty}^\infty
				b( H_b(s\,,x-y)) p_{t-s}(y)\, dy\, ds.
			\]
	\end{enumerate}
\end{theorem}

\begin{remark}
	Before we proceed with a proof,
	we make two remarks:
	\begin{enumerate}
		\item It is possible to adapt the argument of
			Nualart and Pardoux \ycite{NP}*{Proposition 1.6} to deduce that
			the laws of $ H $ and $ H_b$ are
			mutually absolutely continuous with respect to
			one another; see also
			Dalang and Nualart
			\ycite{DN}*{Corollary 5.3} and 
			Dalang, Khoshnevisan, and Nualart
			\ycite{DKN}*{Equations (5.2)
			and (5.3)}. 
			A consequence of this mutual absolute continuity 
			is that $ H_b$ is [H\"older] continuous
			iff $ H $ is.
		\item Theorem \ref{th:heat:semilinear} implies facts that
			are cannot be described by change-of-measure methods.
			For instance, it has
			the striking consequence that with probability one,
			\emph{$ H_b$ and $ H $ blow up in exactly
			the same points}! [This is simply so, because $ H_b
			- H $ is locally bounded.]
			For an example, we mention that the operators
			considered Example
			\ref{example:sharp}, when the parameter $\alpha$ there
			is $\le 2$, lead to discontinuous
			solutions $ H $ that blow up (a.s.) in every open
			subset of $\R_+\times\R$ \cite{MarcusRosen}*{Section 5.3}.
			In those cases, $ H_b$ inherits this property
			as well.
	\end{enumerate}
\end{remark}

\begin{proof}
	We will need the following fact:
	\begin{equation}\label{eq:cont:Hb}
		\text{$ H $ is continuous in probability}.
	\end{equation}
	In fact, we prove that $ H $ is continuous in $L^2(\P)$. 
	
	Owing to \eqref{eq:hEh}, for all $t\ge 0$ and $x,y\in\R$,
	\begin{equation}
		\E\left(\left|  H (t\,,x)- H (t\,,y)
		\right|^2\right) \le \frac{e^{2t}}{2\pi}\int_{-\infty}^\infty
		\frac{1-\cos(\xi|x-y|)}{1+\Re\Psi(\xi)}\, d\xi.
	\end{equation}
	Moreover, according to \eqref{cond:hawkes}, 
	$(1+\Re\Psi)^{-1}\in L^1(\R)$. Therefore,
	the dominated convergence theorem
	implies that for all $T>0$
	and $y\in\R$,
	\begin{equation}\label{eq:leftover1}
		\sup_{t\in[0,T]}\left\|  H (t\,,x)- H (t\,,y)
		\right\|_{L^2(\P)} \to 0 \quad\text{as $x\to y$}.
	\end{equation}
	Similarly, Theorem \ref{th:cont:t} implies that for all $t\ge 0$,
	\begin{equation}
		\lim_{s\to t}
		\sup_{x\in\R}\left\|  H (t\,,x)- H (s\,,x)
		\right\|_{L^2(\P)}^2
		\le \lim_{s\to t}\frac{1}{2\pi}\int_{-\infty}^\infty
		\frac{d\xi}{(1/|t-s|)+\Re\Psi(\xi)},
	\end{equation}
	which is zero by the dominated convergence theorem. This and
	\eqref{eq:leftover1} together imply \eqref{eq:cont:Hb}.
	Now we begin the proof in earnest.
	
	Throughout, we fix
	\begin{equation}
		\text{\rm Lip}_b:= \sup_{x\neq y}\frac{|b(x)-b(y)|}{|x-y|}
		\quad\text{and}\quad
		\lambda:=2\text{\rm Lip}_b.
	\end{equation}
	The condition that $b$ is globally Lipschitz tells
	precisely that $\text{\rm Lip}_b$ and/or $\lambda$ are finite.
	
	Now we begin with a fixed-point scheme: Set $u_0(t\,,x):=0$,
	and define, iteratively for all integers $n\ge 0$,
	\begin{equation}\label{eq:heat:semilinear:picard}
		u_{n+1}(t\,,x) :=  H (t\,,x) + \int_0^t\int_{-\infty}^\infty
		b(u_n(s\,,x-y)) p_{t-s}(y)\, dy\, ds.
	\end{equation}
	
	Consider the processes
	\begin{equation}
		D_{n+1}(x) := \int_0^\infty 
		e^{-\lambda t}\left| u_{n+1}
		(t\,,x)-u_n(t\,,x)
		\right|\,dt\qquad\text{for all $n\ge 0$ 
		and $x\in\R$}.
	\end{equation}
	Also, define $r_\lambda$ to be the $\lambda$-potential density of
	$X$, given by
	\begin{equation}
		r_\lambda(z) := \int_0^\infty p_s(z) e^{-\lambda s}\, ds
		\qquad\text{for all $z\in\R$}.
	\end{equation}
	According to Lemma \ref{lem:heat:tdf}, this is well defined.
	
	The processes $D_1,D_2,\ldots$ satisfy the following recursion:
	\begin{equation}\begin{split}
		D_{n+1}(x) &\le \frac{\text{\rm Lip}_b}{\lambda}(D_n* r_\lambda)(x)\\
		&=\frac12 (D_n* r_\lambda)(x),
	\end{split}\end{equation}
	as can be seen by directly manipulating \eqref{eq:heat:semilinear:picard}.
	We can iterate this to its natural end, and deduce that
	\begin{equation}
		D_{n+1}(x) \le 2^{-n-1}
		(D_0*r_\lambda)(x)\qquad\text{for all $n\ge 0$}.
	\end{equation}
	Since $u_1(t\,,x)-u_0(t\,,x)= H (t\,,x)+tb(0)$,
	\begin{equation}
		D_0(x) \le \int_0^\infty e^{-\lambda t}| H (t\,,x)|
		\, dt + \frac{|b(0)|}{\lambda}.
	\end{equation}
	Thanks to \eqref{eq:cont:Hb}, we can always
	select a measurable modification of
	$(\omega\,,t\,,x)\mapsto H (t\,,x)(\omega)$,
	by the separability theory of Doob
	\ycite{Doob}*{Theorem 2.6, p.\ 61}. Therefore, we can apply the preceding
	to a Lebesgue-measurable modification of $t\mapsto  H (t\,,x)$
	to avoid technical problems \cite{Doob}*{Theorem 2.7, p.\ 62}. 
	In addition, it follows from
	this and convexity that
	\begin{equation}
		\int_{-\infty}^\infty \left\| D_0(x) \right\|_{L^2(\P)}
		e^{-|x|}\, dx \le
		\int_{-\infty}^\infty \int_0^\infty e^{-\lambda t-|x|} 
		\left\|  H (t\,,x)
		\right\|_{L^2(\P)}\, dt\,dx + 
		\frac{2|b(0)|}{\lambda}.
	\end{equation}
	We can apply Proposition \ref{pr:heat:embed}, with its $\lambda$
	replaced by $(4\lambda)^{-1}$ here,
	to find that
	\begin{equation}\begin{split}
		\left\|  H (t\,,x) \right\|_{L^2(\P)}^2
			&\le \frac{e^{\lambda t/2}}{4\pi}\int_{-\infty}^\infty\frac{d\xi}{
			(4\lambda)^{-1}+\Re\Psi(\xi)}\\
		&:= \text{const}
			\cdot e^{\lambda t/2},
	\end{split}\end{equation}
	where the constant depends neither on $t$ nor on $n$.
	Consequently,
	\begin{equation}
		\int_{-\infty}^\infty 
		\|D_0(x)\|_{L^2(\P)}
		e^{-|x|}\,dx <\infty.
	\end{equation}
	Because $\lambda r_\lambda$ is a probability density on $\R$,
	it follows that for all integers $n\ge 0$,
	\begin{equation}
		\int_{-\infty}^\infty
		\left\| D_{n+1}(x) \right\|_{L^2(\P)}
		e^{-|x|}\,dx 
		\le \text{const}\cdot 2^{-n}.
	\end{equation}
	Therefore, $\sum_{n=0}^\infty\int_{-\infty}^\infty
	\|D_n(x)\|_{L^2(\P)}\exp(-|x|)\, dx<\infty$.
	Furthermore, for all $T,k>0$,
	\begin{equation}\begin{split}
		\lim_{n\to\infty} \int_0^T\int_{-k}^k
			\left| u_{n+1}(t\,,x)-u_n(t\,,x)\right|\, dx\, dt
			&\le
			e^{\lambda T+ k}\lim_{n\to\infty} \int_{-k}^k
			D_n(x)e^{-|x|}\, dx\\
		&=0\qquad\text{a.s.}
	\end{split}\end{equation}
	Because $L^1([0\,,T]\times [-k\,,k])$ is complete,
	and since $(\Omega\,,\F,\P)$ is complete,
	standard arguments show
	that there exists a process 
	$u_\infty\in L^1_{\text{\it loc}}(\R_+\times\R)$
	such that
	\begin{equation}
		\lim_{n\to\infty}\int_0^T\int_{-k}^k
		|u_n(t\,,x)-u_\infty(t\,,x)|\,dx\, dt= 0,
	\end{equation}
	almost surely
	for all $T,k>0$.
	Since $b$ is globally Lipschitz, it follows easily from
	this that outside a single set of $\P$-measure zero, 
	\begin{equation}\label{eq:eq:uinfty}
		u_\infty(t\,,x) =  H (t\,,x) + \int_0^t
		\int_{-\infty}^\infty b(u_\infty(s\,,x-y))\, p_{t-s}(y)\, dy\, ds,
	\end{equation}
	simultaneously for almost all $(t\,,x)\in\R_+\times\R$.
	An application of Fubini's theorem
	implies then that the assertion
	\eqref{eq:eq:uinfty} holds a.s.\
	for almost all $(t\,,x)\in\R_+\times\R$.
	[Observe the order of the quantifiers.]
		
	Consider the finite Borel measure 
	\begin{equation}
		\Upsilon(dt\,dx):=e^{-\lambda t-|x|}\,dt\,dx,
	\end{equation}
	defined on $\R_+\times\R$.
	Our proof, so far, contains the fact that
	$u_\infty\in L^1(\Upsilon)$ almost surely. 
	Moreover, thanks to \eqref{eq:eq:uinfty}, if $p\in[1\,,\infty)$  then
	\begin{equation}
		\left\| u_\infty(t\,,x) \right\|_{L^p(\P)}
		\le \left\|  H (t\,,x) \right\|_{L^p(\P)}
		+\sup_z|b(z)|t\qquad
		\text{for all $(t\,,x)\in\R_+\times\R$}.
	\end{equation}
	One of the basic properties of centered Gaussian
	random variables is that their $p$th moment is 
	proportional to their second moment.
	Consequently,
	\begin{equation}
		\left\| u_\infty (t\,,x) \right\|_{L^p(\P)}\le 
		\text{const}\cdot\left( \left\|  H (t\,,x)
		\right\|_{L^2(\P)}+1\right),
	\end{equation}
	which we know to be locally bounded.
	This and the Tonelli theorem together prove that
	$u_\infty\in L^p_{\text{\it loc}}(\R_+\times\R)$ a.s. 
	
	We recall a standard fact from classical analysis:
	\emph{If $f\in L^1(\Upsilon)$,
	then $f$ is continuous in the measure $\Upsilon$}. This means
	that for all $\delta>0$,
	\begin{equation}
		\lim_{\epsilon,\eta\to 0} \Upsilon\left\{ (t\,,x):\
		\big| f(t+\eta\,,x+\epsilon)-f(t\,,x) \big|>\delta\right\}=0,
	\end{equation}
	and follows immediately
	from standard approximation arguments; see 
	the original classic book of Zygmund
	\ycite{Zygmund}*{\S2.201, p.\ 17},
	for instance. Consequently, Lemma \ref{lem:heat:tdf}
	and the already-proved fact that
	$u_\infty\in L^1(\Upsilon)$ a.s.\ together imply
	that $u_\infty$ is
	continuous in measure a.s. Because $b$ is bounded 
	and Lipschitz, the integrability/continuity properties of
	$p_t(x)$, as explained in Lemma \ref{lem:heat:tdf},
	together imply that 
	$(t\,,x)\mapsto \int_0^t\int_{-\infty}^\infty
	b(u_\infty(s\,,x-y)) \, p_{t-s}(y)\,dy\,ds
	=\int_0^t\int_{-\infty}^\infty
	b(u_\infty(s\,,z)) \, p_{t-s}(x-z)\,dz\,ds$
	is a.s.\ continuous, and bounded on $[0\,,T]\times\R$,
	for every nonrandom and fixed $T>0$.
	Now let us \emph{define}
	\begin{equation}
		 H_b(t\,,x) :=  H (t\,,x)
		+\int_0^t\int_{-\infty}^\infty b(u_\infty(s\,,x-y))
		\, p_{t-s}(y)\, dy\, ds.
	\end{equation}
	Thanks to \eqref{eq:eq:uinfty}, 
	\begin{equation}
		\P\left\{  H_b(t\,,x) = u_\infty(t\,,x)
		\right\}=1\qquad\text{for all $(t\,,x)\in\R_+\times\R$.}
	\end{equation}
	Thus, $H_b$ is a modification of $u_\infty$. In
	addition, the Tonelli theorem applies
	to tell us that $ H_b$ inherits the almost-sure
	local integrability property of $u_\infty$.
	That is, $ H_b\in L^p_{\text{\it loc}}
	(\R_+\times\R)$ a.s. Moreover, outside a single null set
	we have
	\begin{equation}
		J(t\,,x)
		=\int_0^t\int_{-\infty}^\infty b(
		u_\infty(s\,,x-y))\, 
		p_{t-s}(y)\, dy\,ds\qquad
		\text{for all $(t\,,x)\in\R_+\times\R$.}
	\end{equation}
	This proves the theorem.
\end{proof}

\end{document}